\documentclass[12pt]{article}
\usepackage{exscale,relsize}
\usepackage{comment}
\usepackage{lscape}
\usepackage{amsmath}
\usepackage[sort,nocompress]{cite}
\usepackage{amsfonts}
\usepackage[hidelinks]{hyperref}
\usepackage{amssymb}
\usepackage{calc}
\usepackage{theorem}
\usepackage{pifont}      
\usepackage{array}
\usepackage{mathtools} 

\usepackage{color}

\usepackage{graphicx} 
\usepackage{float} 
\usepackage{subcaption} 

\newcommand{\Ball}{\mathbb{B}}
\newcommand{\sgn}{\ensuremath{\operatorname{sgn}}}
\newcommand{\N}{\mathbb{N}}   

\newcommand{\dist}{\ensuremath{\operatorname{dist}}} 

\newcommand{\norm}[1]{\left\lVert#1\right\rVert} 
\newcommand{\ip}[2]{\langle#1,#2\rangle} 
\newcommand{\Limsup}{\ensuremath{\operatorname{Limsup}}}

\newcommand{\gph}{\ensuremath{\operatorname{gph}}}

\newcommand{\R}{\ensuremath{\mathbb R}}
\newcommand{\Rn}{\ensuremath{\mathbb R^n}}

\newcommand{\dom}{\ensuremath{\operatorname{dom}}}

\newcommand{\Id}{\ensuremath{\operatorname{Id}}}

\providecommand{\BB}[2]{\mathbb{B}(#1;#2)}

\newcommand{\prox}{\ensuremath{\operatorname{Prox}}}
\newcommand{\proj}{\ensuremath{\operatorname{Proj}}}

\newcommand{\supp}{\ensuremath{\operatorname{supp}}}

\newtheorem{theorem}{Theorem}[section]
\newtheorem{lemma}[theorem]{Lemma}
\newtheorem{fact}[theorem]{Fact}
\newtheorem{corollary}[theorem]{Corollary}
\newtheorem{proposition}[theorem]{Proposition}
\newtheorem{defn}[theorem]{Definition}

\theoremstyle{plain}{\theorembodyfont{\rmfamily}
	\newtheorem{assumption}[theorem]{Assumption}}
\theoremstyle{plain}{\theorembodyfont{\rmfamily}
	}
\theoremstyle{plain}{\theorembodyfont{\rmfamily}
	\newtheorem{algorithm}[theorem]{Algorithm}}
\theoremstyle{plain}{\theorembodyfont{\rmfamily}
	\newtheorem{example}[theorem]{Example}}
\theoremstyle{plain}{\theorembodyfont{\rmfamily}
	\newtheorem{remark}[theorem]{Remark}}
\theoremstyle{plain}{\theorembodyfont{\rmfamily}
	}
\def\proof{\noindent{\it Proof}. \ignorespaces}
\def\endproof{\ensuremath{\quad \hfill \blacksquare}}

\newcommand{\pluss}{{\hskip1pt \raise1pt\vbox{\hrule width6pt \vskip1pt
			\hrule width6pt}\kern-4pt{\lower1pt\hbox{\vrule height6pt \kern1pt\vrule
				height6pt}}\hskip5pt}}

\newcommand{\argmax}{\mathop{\rm argmax}\limits}
\newcommand{\argmin}{\mathop{\rm argmin}\limits}

\begin{document}
\title{A Bregman inertial forward-reflected-backward method for nonconvex minimization}

\author{
         Xianfu\ Wang\thanks{
                Department of Mathematics, Irving K. Barber Faculty of Science,
                University of British Columbia, Kelowna, B.C.\ V1V~1V7, Canada.
                 E-mail: \href{mailto:shawn.wang@ubc.ca}{\texttt{shawn.wang@ubc.ca}}.}\  and
         Ziyuan Wang\thanks{
                 Department of Mathematics, Irving K. Barber Faculty of Science, University of British Columbia, Kelowna, B.C.\ V1V~1V7, Canada.
                 E-mail: \href{mailto:ziyuan.wang@alumni.ubc.ca}{\texttt{ziyuan.wang@alumni.ubc.ca}}.}
                 }

\maketitle

\begin{abstract} \noindent We propose a Bregman inertial forward-reflected-backward~(BiFRB) method for nonconvex composite problems. Our analysis relies on a novel approach that imposes general conditions on implicit merit function parameters, which yields a stepsize condition that is independent of inertial parameters. In turn,
a question of Malitsky and Tam regarding whether FRB can be equipped with a Nesterov-type acceleration is resolved. Assuming the generalized concave
Kurdyka-\L ojasiewicz property of a quadratic regularization of the objective, we obtain sequential convergence of BiFRB, as well as convergence rates on
both the function value and actual sequence. We also present formulae for the Bregman subproblem, supplementing not only BiFRB but also the work of Bo{\c{t}}-Csetnek-L{\'a}szl{\'o} and Bo{\c{t}}-Csetnek. Numerical simulations are conducted to evaluate the performance of our proposed algorithm.
\end{abstract}


\section{Introduction}
Consider the inclusion problem of maximally monotone operators:
\begin{equation}\label{inclusion problem}
	\text{find }x\in\Rn\text{ such that }0\in (A+B)(x),
\end{equation}
where $A:\Rn\to 2^{\Rn}$ and $B:\Rn\to\Rn$ are (maximally) monotone operators with $B$ being Lipschitz continuous. Behavior of splitting methods for solving~(\ref{inclusion problem}) is well-understood; see, e.g., the recent forward-reflected-backward method by Malitsky and Tam~\cite{malitsky2020FRB}, and~\cite[Chapters 26, 28]{BC}, which includes classical algorithms such as the forward-backward method~\cite{combettes2005signal} with $B$ being cocoercive, the Douglas-Rachford method~\cite{lions1979splitting,douglas1956numerical} and Tseng's method~\cite{tseng2000modified}.

\emph{The goal of this paper is to propose a Bregman inertial forward-reflected-backward method~(Algorithm~\ref{BiFRB}) for solving the nonconvex composite problem}
 \begin{equation}\label{optimization prblem}
 	\min_{x\in\Rn} F(x)=f(x)+g(x),		
 \end{equation}
 where $f:\Rn\to\overline{\R}=(-\infty,\infty]$ is proper lower semicontinuous (lsc),
 and $g:\Rn\to\R$ has a Lipschitz continuous gradient, which corresponds to~(\ref{inclusion problem}) with $A=\partial f$ and $B=\nabla g$ when $f$ and $g$ are both convex. Informally, the proposed algorithm computes
 \begin{equation}\label{formula: informal algorithm}
 	x_{k+1}\in (\nabla h+\lambda_k \partial f)^{-1}(x_k-2\lambda_k  \nabla g(x_k)+\lambda_k  \nabla g(x_{k-1})+\alpha_k(x_k-x_{k-1})),
 \end{equation}
 for some stepsize $\lambda_k>0$, inertial parameter $\alpha_k\in[0,1)$ and kernel $h:\Rn\to\R$. In particular,~(\ref{formula: informal algorithm}) reduces to the inertial forward-reflected-backward method studied in~\cite[Corollary 4.4]{malitsky2020FRB} when $h=\norm{\cdot}^2/2$, $f,g$ are convex and inertial parameter is fixed.
 Before stating our main contribution, let us note that several splitting algorithms have been extended to nonconvex setting for solving~(\ref{optimization prblem}); see, e.g., the Douglas-Rachford method~\cite{li2016douglas}, forward-backward method~\cite{attouch2013convergence,bolte2018first}, inertial forward-backward methods~\cite{ochs2014ipiano,boct2016inertiala}, the Peaceman-Rachford method~\cite{li2017peaceman}, inertial Tseng's method~\cite{boct2016inertial}, and forward-reflected-backward method~\cite{wang2021malitsky}.

Let us summarize the \emph{main contributions} of this paper. Our convergence analysis relies on the
 generalized concave Kurdyka-\L ojasiewicz property~(Definition~\ref{def: g-concave KL}) and {a novel framework for analyzing a sequence of implicit merit functions~(Definition~\ref{def: merit}) through a general condition~(Assumption~\ref{assmption: p_k and M_1,k}) on merit function parameters}. In turn, we derive {global sequential convergence to a stationary point of $F$ with
 an explicit stepsize condition that is independent of inertial parameters}, whereas stepsize assumptions in current literature are tangled with inertial parameters
 or even implicit; see Remark~\ref{rem: parameter condition advantage}(i). Notably, {such an independence result resolves a question of Malitsky and
 Tam~\cite[Section 7]{malitsky2020FRB} regarding whether FRB can be adapted to incorporate a Nesterov-type acceleration; see Remark~\ref{rem: parameter condition advantage}(ii).}
 {Convergence rate analysis on both the function value and actual sequence} is also carried out. Moreover, we
 {provide several formulae of the associated Bregman subproblem}, which to the best of our knowledge are the first of its kind among present results
 devoting to Bregman extensions of splitting algorithms with the same kernel; see, e.g.,~\cite{boct2016inertial,boct2016inertiala}.

The paper is organized as follows.
We begin with notation and background in Section~\ref{sec:Notation and preliminaries}. The proposed algorithm and its abstract convergence are studied in Section~\ref{sec:Bregman inertial FRB splitting method}. Then, in Section~\ref{sec:Stepsize strategies }, we present suitable stepsize rules under which the abstract convergence holds. Function value and sequential convergence rate results are discussed in Section~\ref{sec:rate}, followed by Bregman subproblem formulae in Section~\ref{sec:subproblem formulae}. Numerical simulations are conducted in Section~\ref{sec:numerical}. We end this paper in Section~\ref{sec:conclusion and future work} with remarks on our contribution and possible directions for future research.

 \section{Notation and preliminaries}\label{sec:Notation and preliminaries}
 Throughout this paper, $\Rn$ is the standard Euclidean space equipped with inner product~$\ip{x}{y}=x^Ty$ and the Euclidean norm~$\norm{x}=\sqrt{\ip{x}{x}}$ for $x, y\in\Rn$.
 Let~$\N=\{0,1,2,\ldots\}$ and~$\N^*=\{-1\}\cup\N$. The open ball centered at $\bar{x}$ with radius $r$ is denoted by~$\BB{\bar{x}}{r}$.
 The distance function of a subset $K\subseteq \Rn$ is $\dist(\cdot,K):\Rn\rightarrow\overline{\R}=(-\infty,\infty]$,
 \[x\mapsto\dist(x,K)=\inf\{\norm{x-y}:y\in K\},\]
 where $\dist(x,K)\equiv\infty$ if $K=\emptyset$. For $f:\Rn\to\overline{\R}$ and $r_1,r_2\in[-\infty,\infty]$, we set $[r_1<f<r_2]=\{x\in\Rn:r_1<f(x)<r_2\}$. We say a function $f:\Rn\to\overline{\R}$ is coercive if $\lim_{\norm{x}\to\infty}f(x)=\infty$. A proper, lsc function $f:\Rn\to\overline{\R}$ is prox-bounded if there exists $\lambda>0$ such that $f+\norm{\cdot}^2/(2\lambda)$ is bounded below; see, e.g.,~\cite[Exercise 1.24]{rockwets}.
 The supremum of the set of all such $\lambda$ is the threshold $\lambda_f$ of prox-boundedness for $f$. The indicator function of a set $K\subseteq\Rn$ is $\delta_K(x)=0$ if $x\in K$ and $\delta_K(x)=\infty$ otherwise.

 We will use frequently the following concepts from variational analysis; see, e.g., \cite{rockwets, mor2005variational}.

 \begin{defn}\label{Defn:limiting subdifferential}
 	Let $f:\Rn\rightarrow\overline{\R}$ be a proper function and $K\subseteq\Rn$ a nonempty closed set. We say that

(i) $v\in\Rn$ is a \textit{Fr\'echet subgradient} of $f$ at $\bar{x}\in\dom f$, denoted by $v\in\hat{\partial}f(\bar{x})$, if
 		\begin{equation}\label{Formula:frechet subgradient inequality}
 			f(x)\geq f(\bar{x})+\ip{v}{x-\bar{x}}+o(\norm{x-\bar{x}}).
 		\end{equation}

(ii) $v\in\Rn$ is a \textit{limiting subgradient} of $f$ at $\bar{x}\in\dom f$, denoted by $v\in\partial f(\bar{x})$, if
 		\begin{equation}\label{Formula:limiting subgraident definition}
 			v\in\{v\in\Rn:\exists x_k\xrightarrow[]{f}\bar{x},\exists v_k\in\hat{\partial}f(x_k),v_k\rightarrow v\},
 		\end{equation}
 		where $x_k\xrightarrow[]{f}\bar{x}\Leftrightarrow x_k\rightarrow\bar{x}\text{ and }f(x_k)\rightarrow f(\bar{x})$. Moreover, we set $\dom\partial f=\{x\in\Rn:\partial f(x)\neq\emptyset\}$. We say that $\bar{x}\in\dom\partial f$ is a stationary point, if $0\in\partial f(\bar{x})$.
 		
(iii) The Fr\'echet and limiting normal cones to $K$ at $\bar{x}$ are $\hat{N}_K(\bar{x})=\hat{\partial}\delta_K(\bar{x})$ and $N_K(\bar{x})=\partial\delta_K(\bar{x})$, respectively, if $\bar{x}\in K$; otherwise $\hat{N}_K(\bar{x})=N_K(\bar{x})=\emptyset$.
 \end{defn}
If $K\subseteq\Rn$ is a nonempty convex set, then $N_K(\bar{x})=\hat{N}_K(\bar{x})=\{v\in\Rn:(\forall x\in K)~\ip{v}{x-\bar{x}}\leq0\}$; see, e.g.,~\cite{rockwets}.

The proximal mapping of a
 proper and lsc $f:\Rn\to\overline{\R}$ with parameter $\lambda>0$ is
 \[ (\forall x\in\Rn)\ \prox_{\lambda f}(x)=\argmin_{y\in\Rn}\Big\{f(y)+\frac{1}{2\lambda}\norm{x-y}^2\Big\},\]
 which is the resolvent $J_{\lambda A}$ with $A=\partial f$ when $f$ is convex; see, e.g.,~\cite{BC}. The Bregman distance induced by a differentiable $h:\Rn\to\R$ is \[D_h: \Rn\times\Rn\to\R: (x,y)\mapsto h(x)-h(y)-\ip{x-y}{\nabla h(y)}.\]

For $\eta\in(0,\infty]$, denote by $\Phi_\eta$ the class of functions $\varphi:[0,\eta)\rightarrow\R_+$ satisfying the following conditions: (i) $\varphi(t)$ is right-continuous at $t=0$ with $\varphi(0)=0$; (ii) $\varphi$ is strictly increasing on $[0,\eta)$. The following concept will be the key in our convergence analysis.

\begin{defn}\emph{\cite{wang2020}}\label{def: g-concave KL} Let $f:\Rn\rightarrow\overline{\mathbb{R}}$
	be proper and lsc.  Let $\bar{x}\in\dom\partial f$ and $\mu\in\R$, and
	let $V\subseteq\dom\partial f$ be a nonempty subset.
	
	(i) We say that $f$ has the pointwise generalized concave Kurdyka-\L ojasiewicz~(KL) property at $\bar{x}\in\dom\partial f$, if there exist neighborhood $U\ni\bar{x}$, $\eta\in(0,\infty]$ and concave $\varphi\in\Phi_\eta$, such that for all $x\in U\cap[0<f-f(\bar{x})<\eta]$,
	\begin{equation}\label{g-concave KL inequality}
		\varphi^\prime_-\big(f(x)-f(\bar{x})\big)\cdot\dist\big(0,\partial f(x)\big)\geq1,
	\end{equation}
	where $\varphi_-^\prime$ denotes the left derivative. Moreover, $f$ is a generalized concave KL function if it has the generalized concave KL property at every $x\in\dom\partial f$.
	
	(ii) Suppose that $f(x)=\mu$ on $V$. We say $f$ has the setwise\footnote{We shall omit adjectives ``pointwise" and ``setwise" whenever there is no ambiguity.} generalized concave Kurdyka-\L ojasiewicz property on $V$, if there exist $U\supset V$, $\eta\in(0,\infty]$ and concave $\varphi\in\Phi_\eta$ such that for every $x\in U\cap[0<f-\mu<\eta]$,
	\begin{equation}\label{uniform g-concave KL inequality}
		\varphi^\prime_-\big(f(x)-\mu\big)\cdot\dist\big(0,\partial f(x)\big)\geq1.
	\end{equation}
\end{defn}

Generalized concave KL functions are ubiquitous. The class of semialgebraic functions is one particularly rich resource. We refer the interested reader to~\cite{bolte2010survey,Attouch2010,bolte2007} for more examples; see also the fundamental work of \L ojasiewicz~\cite{Lojas1963} and Kurdyka~\cite{Kur98}.

 \begin{defn} (i) A set $E\subseteq\mathbb{R}^n$ is called semialgebraic if there exist finitely many polynomials $g_{ij}, h_{ij}:\mathbb{R}^n\rightarrow\mathbb{R}$ such that	
 	$E=\bigcup_{j=1}^p\bigcap_{i=1}^q\{x\in\mathbb{			R}^n:g_{ij}(x)=0\text{ and }h_{ij}(x)<0\}.	$
 	
 	(ii) A function $f:\mathbb{R}^n\rightarrow\overline{\R}$ is called semialgebraic if its graph $\gph f=\{(x,y)\in\mathbb{R}^{n+1}:f(x)=y\}$ is semialgebraic.
 \end{defn}

 \begin{fact}\label{Fact:semi-algebraic functions are KL} \emph{\cite[Corollary 16]{bolte2007}} Let $f:\mathbb{R}^n\rightarrow\overline{\R}$ be a proper and lsc function and let $\bar{x}\in\dom\partial f$. If $f$ is semialgebraic, then it has the concave KL property at $\bar{x}$ with $\varphi(t)=c\cdot t^{1-\theta}$ for some $c>0$ and $\theta\in(0,1)$.
 \end{fact}

\section{The BiFRB method and its abstract convergence}\label{sec:Bregman inertial FRB splitting method}
\subsection{The BiFRB method}
In the remainder of this paper, suppose that the following hold:
\begin{assumption}\label{Standing assumption}
	
	(i) $f:\Rn\to\overline{\R}$ is proper, lsc and  prox-bounded with threshold $\lambda_f>0$.
	
	(ii) $g:\Rn\to\R$ has Lipschitz continuous gradient with constant $L_{\nabla g}$.
	
	(iii) $h:\Rn\to\R$ is $\sigma$-strongly convex and has Lipschitz continuous gradient with constant $L_{\nabla h}$.
	
	(iv) The objective $F=f+g$ is bounded below.
\end{assumption}
We now propose the Bregman inertial forward-reflected-backward (BiFRB) algorithm.

\begin{algorithm}[BiFRB]\label{BiFRB}~~

1. Initialization: Pick $x_{-1}, x_0\in\Rn$. Let $0<\underline{\lambda}\leq\overline{\lambda}\leq\lambda_f$.  Choose $\underline{\lambda}\leq\lambda_{-1}\leq\overline{\lambda}$.

	2. For $k\in\N$, choose $\underline{\lambda}\leq\lambda_{k}\leq\overline{\lambda}$ and $0\leq\alpha_{k}<1$. Compute
\begin{align}
					&y_k=x_k+\lambda_{k-1}\big(\nabla g(x_{k-1})-\nabla g(x_k) \big),\label{Algorithm: y_k}\\
					&x_{k+1}\in\argmin_{x\in\Rn}\bigg\{f(x)+
					\ip{x-y_k}{\nabla g(x_k)+\frac{\alpha_k}{\lambda_k}(x_{k-1}-x_k)}+\frac{1}{\lambda_k}D_h(x,y_k)\bigg\}.\label{Algorithm: x_{k+1}}
		\end{align}

\end{algorithm}

The Bregman distance $D_h$ is proportional to square of the Euclidean distance.
 \begin{lemma}\label{lem: Bregman distance}\emph{\cite[Exercise 17.5, Theorem 18.15]{BC}} Let $h:\Rn\to\R$ be $\sigma$-strongly convex and has Lipschitz continuous gradient with constant $L_{\nabla h}$. Then \[(\forall x,y\in\Rn)~\frac{\sigma}{2}\norm{x-y}^2\leq D_h(x,y)\leq\frac{L_{\nabla h}}{2}\norm{x-y}^2.\]
\end{lemma}


Assumption~\ref{Standing assumption}(iii) has been used in the literature to obtain non-Euclidean extension of splitting algorithms; see, e.g.,~\cite{boct2016inertial,boct2016inertiala}.  Recall that the Fenchel conjugate of $f:\Rn\to\overline{\R}$ is $f^*=\sup_{y\in\Rn}\{\ip{\cdot}{y}-f(y)\}$ and the Moreau envelope of $f:\Rn\to\overline{\R}$ with parameter $\lambda>0$ is \[(\forall x\in\Rn)~M_f^\lambda(x)=\inf_{y\in\Rn}\bigg\{f(y)+\frac{1}{2\lambda}\norm{x-y}^2\bigg\}.\]
Below is a systemic way to construct kernel $h$ via Fenchel conjugate and Moreau envelope.

\begin{proposition}\label{prop: kernel construction} Let $L,\sigma>0$ and let $p:\Rn\to\overline{\R}$ be proper, lsc and convex. Define $q=\norm{\cdot}^2/2$. Then the following hold:
	
(i)\emph{~\cite[Corollary 18.18]{BC}} $p$ has $L$-Lipschitz gradient if and only if $p$ is the Moreau envelope with parameter $\lambda=1/L$ of a proper lsc and convex function $(p^*-\lambda q)^*$. Therefore, in particular, if $p^*$ is $1/L$-strongly convex then $p$ has $L$-Lipschitz gradient.

(ii) If $p$ has $L$-Lipschitz gradient, then $h=p+\sigma q$ is $\sigma$-strongly convex with $(L+\sigma)$-Lipschitz gradient.
\end{proposition}

\begin{example} Let $h(x)=\sqrt{1+\norm{x}^2}+\norm{x}^2/2$. Then $h$ satisfies Assumption~\ref{Standing assumption}(iii) with $\sigma=1$ and $L_{\nabla h}=2$.
\end{example}
\proof Define $(\forall x\in\Rn)~p(x)=\sqrt{1+\norm{x}^2}$. Then $p^*$ is $1$-strongly convex; see, e.g.,~\cite[Section 4.4, Example 5.29]{beck2017first}. Applying Proposition~\ref{prop: kernel construction} completes the proof.\endproof

We also assume the following throughout the paper, under which the proposed algorithm is well-defined.
\begin{assumption}\label{assumption: well-defined} Suppose that one of the following holds:

(i) The kernel $h$ has strong convexity modulus $\sigma\geq1$.

(ii) Function $f$ has prox-threshold $\lambda_f=\infty$.
\end{assumption}
\begin{remark}The threshold $\lambda_f=\infty$ when $f$ is convex or bounded below.
\end{remark}
\begin{proposition} Let $u,v\in\Rn$ and suppose that Assumption~\ref{assumption: well-defined} holds.Then for $0<\lambda<\lambda_f$, the function $x\mapsto f(x)+\ip{x-u}{v}+D_h(x,u)/\lambda$ is coercive. Consequently, the set $$\argmin_{x\in\Rn} \bigg\{f(x)+\ip{x-u}{v}+\frac{1}{\lambda}D_h(x,u)\bigg\}\neq\emptyset.$$
\end{proposition}
\proof We claim that $x\mapsto f(x)+\ip{x}{\omega}$ is prox-bounded with $\lambda_f$ for an arbitrary $\omega\in\Rn$. Indeed,  \[\liminf_{\norm{x}\to\infty}\frac{f(x)+\ip{x}{\omega}}{\norm{x}^2}=\liminf_{\norm{x}\to\infty}\frac{f(x)}{\norm{x}^2}+\lim_{\norm{x}\to\infty}\frac{\ip{x}{\omega}}{\norm{x}^2} =\liminf_{\norm{x}\to\infty}\frac{f(x)}{\norm{x}^2}.\]
Thus invoking~\cite[Excercise 1.24]{rockwets} proves the claim.

By our claim, $x\mapsto f(x)+\ip{x}{v-\sigma u/\lambda}$ is prox-bounded with threshold $\lambda_f$. Let $\lambda^\prime$ be such that $\lambda<\lambda^\prime<\lambda_f$. Then $\lambda^\prime/\sigma<\lambda_f$ under Assumption~\ref{assumption: well-defined}. Therefore
\begin{align*}
&f(x)+\ip{x-u}{v}+\frac{1}{\lambda}D_h(x,u)\geq f(x)+\ip{x-u}{v}+\frac{\sigma}{2\lambda}\norm{x-u}^2-\frac{\sigma}{2\lambda^\prime}\norm{x}^2+\frac{\sigma}{2\lambda^\prime}\norm{x}^2\\
&\geq\frac{\sigma}{2}\left(\frac{1}{\lambda}-\frac{1}{\lambda^\prime}\right)\norm{x}^2+f(x)+\ip{x}{v-\frac{\sigma}{\lambda}u}+\frac{\sigma}{2\lambda^\prime}\norm{x}^2-\ip{u}{v}\\
&\geq \frac{\sigma}{2}\left(\frac{1}{\lambda}-\frac{1}{\lambda^\prime}\right)\norm{x}^2+\inf_{x\in\Rn}\bigg\{f(x)+\ip{x}{v-\frac{\sigma}{\lambda}u}+\frac{\sigma}{2\lambda^\prime}\norm{x}^2\bigg\}-\ip{u}{v} \to\infty,
\end{align*}
as $\norm{x}\to\infty$, where the infimum in the last inequality is finite thanks to prox-boundedness of $x\mapsto f(x)+\ip{x}{v-\sigma u/\lambda}$, which completes the proof.\endproof

Setting the kernel $h=\norm{\cdot}/2$ gives Algorithm~\ref{iFRB}, which is an inertial forward-reflected-backward method~(iFRB).

\begin{algorithm}[iFRB]~\label{iFRB}~~
	
	1. Initialization: Pick $x_{-1}, x_0\in\Rn$. Let $0<\underline{\lambda}\leq\overline{\lambda}$.  Choose $\underline{\lambda}\leq\lambda_{-1}\leq\overline{\lambda}$.
	
	2. For $k\in\N$, choose $\underline{\lambda}\leq\lambda_{k}\leq\overline{\lambda}$ and $0\leq\alpha_{k}<1$. Compute
	\begin{align*}
		&y_k=x_k+\lambda_{k-1}\big(\nabla g(x_{k-1})-\nabla g(x_k) \big),\\
		&x_{k+1}\in\argmin_{x\in\Rn}\bigg\{f(x)+\ip{x-y_k}{\nabla g(x_k)+\frac{\alpha_k}{\lambda_k}(x_{k-1}-x_k)}+\frac{1}{2\lambda_k}\norm{x-y_k}^2\bigg\}.
	\end{align*}
\end{algorithm}
Convergence of iFRB with fixed stepsize and inertial parameter was studied
in~\cite{malitsky2020FRB} in the presence of convexity. However, its behavior in nonconvex setting is still not understood, which will be covered in our analysis.

\subsection{BiFRB merit function and its parameter rules}
In the remainder of this paper, unless specified otherwise, let $p_{-1}\in\R$ and define for $k\in\N$
\begin{align}
	p_{k}&=\left(\frac{(L_{\nabla h}-\sigma)L^2_{\nabla g}}{2}\right)\frac{\lambda^2_{k-1}}{\lambda_k}+\frac{1}{2}\left(\frac{\sigma}{\lambda_k}-L_{\nabla g} \right)-p_{k-1},\label{def: p_k}\\
	M_{1,k}&=p_{k-1}-\frac{\alpha_k+\sigma L_{\nabla g}\lambda_{k-1}}{2\lambda_k}-\frac{(L_{\nabla h}-\sigma)L^2_{\nabla g}\lambda_{k-1}^2}{2\lambda_k},\label{def:M_1,k}
\end{align}
where $(\lambda_k)_{k\in\N^*}$ is a stepsize sequence of BiFRB; $(\alpha_k)_{k\in\N}$ is a sequence of inertial parameters; and $L_{\nabla h}, L_{\nabla g}, \sigma$ are parameters specified in Assumption~\ref{Standing assumption}. The following concept plays a central role in our analysis.
\begin{defn}\label{def: merit} Let $p\in\R$. Define $H_p:\Rn\times\Rn\to\overline{\R}$ by $$(\forall (x,y)\in\Rn\times\Rn)~H_p(x,y)=F(x)+p\norm{x-y}^2.$$ Then we say that $H_p$ is a quadratic merit function with parameter $p$. Let $(p_k)_{k\in\N^*}$ be given by~(\ref{def: p_k}).  Then $H_{p_k}$ is the $k^{th}$ BiFRB merit function. If $p_k\leq\bar{p}$ for some $\bar{p}>0$, then $H_{\bar{p}}$ is a dominating BiFRB merit function.
\end{defn}
\begin{remark} Quadratic merit functions of various forms appear frequently in literature that devotes to extending splitting methods to nonconvex setting; see, e.g.,~\cite{boct2016inertial,boct2016inertiala,wang2021malitsky,ochs2014ipiano,ochs19unifying} and the references therein.
\end{remark}

Despite the recursion~(\ref{def: p_k}), we shall see soon that $(p_k)_{k\in\N^*}$ and $(M_{1,k})_{k\in\N^*}$ are indeed implicit in our analysis. Imposing the following novel condition on $(p_k)_{k\in\N^*}$ and $(M_{1,k})_{k\in\N^*}$, one gets abstract convergence of BiFRB, which in turn yields less restrictive stepsize rules; see Section~\ref{sec:Stepsize strategies } and in particular Remark~\ref{rem: parameter condition advantage}.
\begin{assumption}[general rules of BiFRB merit function parameter]\label{assmption: p_k and M_1,k} Let $(p_k)_{k\in\N^*}$ and $(M_{1,k})_{k\in\N}$ be given as in~(\ref{def: p_k}) and~(\ref{def:M_1,k}). Suppose that the following hold:
	
	(i) $\liminf_{k\to\infty}p_k\geq0$ and $\liminf_{k\to\infty} M_{1,k}>0$.
	
	(ii) There exists $\bar{p}$ such that $p_k\leq \bar{p}$ for all $k$.
\end{assumption}

\begin{remark} Assumption~\ref{assmption: p_k and M_1,k} generalizes quadratic merit function properties in the present literature; see,~e.g.,~\cite[Lemma 3.2]{boct2016inertial},~\cite[Lemma 6]{boct2016inertiala},~\cite[Assumption (H1)]{ochs19unifying} and~\cite[Algorithm 5, Lemma 4.6, Proposition 4.7]{ochs2014ipiano}.
\end{remark}

\subsection{Abstract function value convergence}\label{sec: basic properties}
In this section, we study function value convergence of BiFRB under Assumption~\ref{assmption: p_k and M_1,k}; see Section~\ref{sec:Stepsize strategies } for concrete stepsize rules under which Assumption~\ref{assmption: p_k and M_1,k} holds.
\begin{lemma}\label{lem: merit function descent}  Let $(x_k)_{k\in\N^*}$ be a sequence generated by BiFRB and define $z_k=(x_{k+1},x_k)$ for $k\in\N^*$.  Then the following hold:

(i) For all $k\in\N$,
\begin{equation}\label{fom: descent property}
H_{p_{k}}(z_{k})\leq H_{p_{k-1}}(z_{k-1})-M_{1,k}\norm{z_k-z_{k-1}}^2.
\end{equation}

(ii) Suppose that Assumption~\ref{assmption: p_k and M_1,k}(i) holds. Then the sequence $\left(H_{p_k}(z_k)\right)_{k\in\N}$ is decreasing and $\sum_{k=0}^\infty M_{1,k}\norm{z_k-z_{k-1}}^2<\infty$. Consequently, $\left(H_{p_k}(z_k)\right)_{k\in\N}$ is convergent and $\lim_{k\to\infty}\norm{z_k-z_{k-1}}=0$.
\end{lemma}
\proof
(i) By using the iterative scheme~(\ref{Algorithm: x_{k+1}}),
\begin{align}
f(x_{k+1})&\leq f(x_k)+\ip{x_k-x_{k+1}}{\nabla g(x_k)}+\frac{\alpha_k}{\lambda_k}\ip{x_k-x_{k+1}}{x_{k-1}-x_k}\nonumber \\
&+\frac{1}{\lambda_k}\left(D_h(x_k,y_k)-D_h(x_{k+1},y_k)\right)\nonumber\\
&\leq f(x_k)+\ip{x_k-x_{k+1}}{\nabla g(x_k)}+\frac{\alpha_k}{\lambda_k}\ip{x_k-x_{k+1}}{x_{k-1}-x_k}\nonumber \\
&+\frac{L_{\nabla h}}{2\lambda_k}\norm{x_k-y_k}^2-\frac{\sigma}{2\lambda_k}\norm{x_{k+1}-y_k}^2,\label{fom: f descent}
\end{align}
where the last inequality is implied by Lemma~\ref{lem: Bregman distance}.  Invoking the descent lemma~\cite[Lemma 2.64]{BC} to $g$ yields
\begin{equation}\label{fom: g descent lemma}
g(x_{k+1})\leq g(x_k)+\ip{x_{k+1}-x_k}{\nabla g(x_k)}+\frac{L_{\nabla g}}{2}\norm{x_{k+1}-x_k}^2.
\end{equation}
Adding~(\ref{fom: f descent}) and~(\ref{fom: g descent lemma}), one has
\begin{align}
F(x_{k+1})&\leq F(x_k)+\frac{\alpha_k}{\lambda_k}\ip{x_k-x_{k+1} }{x_{k-1}-x_k}+\frac{L_{\nabla h}}{2\lambda_k}\norm{x_k-y_k}^2\nonumber \\
&-\frac{\sigma}{2\lambda_k}\norm{x_{k+1}-y_k}^2+\frac{L_{\nabla g}}{2}\norm{x_{k+1}-x_k}^2\nonumber\\
&\leq F(x_k)+\left(\frac{\alpha_k}{2\lambda_k}+\frac{L_{\nabla g }}{2}\right)\norm{x_{k+1}-x_k}^2+\frac{\alpha_{k}}{2\lambda_k}\norm{x_{k-1}-x_k}^2\nonumber \\
&+\frac{L_{\nabla h}}{2\lambda_k}\norm{x_k-y_k}^2-\frac{\sigma}{2\lambda_k}\norm{x_{k+1}-y_k}^2.\label{fom:xxx}
\end{align}
Next we estimate
\begin{align*}
& \frac{L_{\nabla h}}{2\lambda_k}\norm{x_k-y_k}^2-\frac{\sigma}{2\lambda_k}\norm{x_{k+1}-y_k}^2
=\frac{L_{\nabla h}-\sigma}{2\lambda_k}\norm{x_k-y_k}^2-\frac{\sigma}{2\lambda_k}\left(\norm{x_{k+1}-y_k}^2-\norm{x_k-y_k}^2\right)\\
&=\frac{(L_{\nabla h}-\sigma)\lambda_{k-1}^2}{2\lambda_k}\norm{\nabla g(x_{k})-\nabla g(x_{k-1}) }^2-\frac{\sigma}{2\lambda_k}\norm{x_{k+1}-x_k}^2\\
&-\frac{\sigma\lambda_{k-1}}{\lambda_k}\ip{x_{k+1}-x_k}{\nabla g(x_{k})-\nabla g(x_{k-1})}\\
&\leq \frac{(L_{\nabla h}-\sigma)L_{\nabla g}^2\lambda_{k-1}^2}{2\lambda_k}\norm{x_k-x_{k-1}}^2+\frac{\sigma L_{\nabla g}\lambda_{k-1}}{\lambda_k}\norm{x_{k+1}-x_k}\norm{x_k-x_{k-1}}-\frac{\sigma}{2\lambda_k}\norm{x_{k+1}-x_k}^2\\
&\leq\frac{(L_{\nabla h}-\sigma)L_{\nabla g}^2\lambda_{k-1}^2+\sigma L_{\nabla g}\lambda_{k-1}}{2\lambda_k}\norm{x_k-x_{k-1}}^2+\frac{\sigma(L_{\nabla g}\lambda_{k-1}-1)}{2\lambda_k}\norm{x_{k+1}-x_k}^2.
\end{align*}
Thus, inequality~(\ref{fom:xxx}) further implies that
\begin{align}\label{formula: fdd}
F(x_{k+1})&\leq F(x_k)+\left(\frac{\sigma(L_{\nabla g}\lambda_{k-1}-1)+\alpha_k}{2\lambda_k}+\frac{L_{\nabla g}}{2} \right)\norm{x_{k+1}-x_k}^2\nonumber \\
&+\frac{(L_{\nabla h}-\sigma)L_{\nabla g}^2\lambda_{k-1}^2+\sigma L_{\nabla g}\lambda_{k-1}+\alpha_k}{2\lambda_k}\norm{x_k-x_{k-1}}^2.
\end{align}
Adding $p_{k}\norm{x_{k+1}-x_k}^2$ to the above inequality gives
\begin{align*}
H_{p_{k}}(z_{k})&=F(x_{k+1})+p_{k}\norm{x_{k+1}-x_k}^2\\
&\leq F(x_k)+\left(\frac{\sigma(L_{\nabla g}\lambda_{k-1}-1)+\alpha_k}{2\lambda_k}+\frac{L_{\nabla g}}{2}+p_{k}\right)(\norm{z_k-z_{k-1}}^2-\norm{x_k-x_{k-1}}^2)\\
&+\frac{(L_{\nabla h}-\sigma)L_{\nabla g}^2\lambda_{k-1}^2+\sigma L_{\nabla g}\lambda_{k-1}+\alpha_k}{2\lambda_k}\norm{x_k-x_{k-1}}^2\\
&= F(x_k)+p_{k-1}\norm{x_k-x_{k-1}}^2-M_{1,k}\norm{z_k-z_{k-1}}^2\\
&=H_{p_{k-1}}(z_{k-1})-M_{1,k}\norm{z_k-z_{k-1}}^2,
\end{align*}
where the first inequality follows from the identity $\norm{z_k-z_{k-1}}^2=\norm{x_{k+1}-x_k}^2+\norm{x_k-x_{k-1}}^2$.

(ii) Assume without loss of generality that $p_k\geq0$ and $M_{1,k}>0$ for $k\in\N$. Thus by statement(i) and Assumption~\ref{Standing assumption}, the sequence $\left(H_{p_k}(z_k)\right)_{k\in\N}$ is decreasing and bounded below, therefore it is convergent. Summing inequality~(\ref{fom: descent property}) from $k=0$ to $k=n-1$, we have
\begin{align*}
\sum_{k=0}^{n-1}M_{1,k}\norm{z_k-z_{k-1}}^2&\leq \sum_{k=0}^{n-1}\big(H_{p_k}(z_k)-H_{p_{k+1}}(z_{k+1})\big)= H_{p_0}(z_0)-H_{p_n}(z_n)\\
&\leq H_{p_0}(z_0)-F(x_{n+1})\leq H_{p_0}(z_0)-\inf_{x\in\Rn}F(x)<\infty,
\end{align*}
which together with the assumption that $\liminf_{k\to\infty} M_{1,k}>0$ implies that $\sum_{k=0}^\infty\norm{z_k-z_{k-1}}^2<\infty$.\endproof

The following lemma is an immediate consequence of subdifferential sum rules~\cite[Exercise 8.8, Proposition 10.5]{rockwets}.
\begin{lemma}\label{lem: subdifferential H} Let $p\in\R$. Then for every $(x,y)\in\Rn\times\Rn$,
\[\partial H_p(x,y)=\{\partial f(x)+\nabla g(x)+2p(x-y)\}\times \{2p(y-x)\}. \]
\end{lemma}

The lemma below provides a lower bound for consecutive iterates gap.

\begin{lemma}\label{lem: subdifferential bound} Let $(x_k)_{k\in\N^*}$ be a sequence generated by BiFRB, $z_k=(x_{k+1},x_k)$ for $k\in\N^*$, and $p\geq0$.  For $k\in\N$, define $u_{k}=(\nabla h(y_k)-\nabla h(x_{k+1}))/\lambda_k-\nabla g(x_k)+\alpha_k(x_k-x_{k-1})/\lambda_k$,
\[A_k=u_{k}+\nabla g(x_{k+1})+2p(x_{k+1}-x_k), B_k=2p(x_k-x_{k+1}). \]
Set \[M_{2,k}=\sqrt{2}\max\bigg\{\frac{L_{\nabla h}}{\lambda_k}+L_{\nabla g}+6p,  \frac{L_{\nabla h}L_{\nabla g}\lambda_{k-1}+1}{\lambda_k}\bigg\},\]
and
\[M_{2}=\sqrt{2}\max\bigg\{\frac{L_{\nabla h}}{\underline{\lambda}}+L_{\nabla g}+6p,  \frac{L_{\nabla h}L_{\nabla g}\overline{\lambda}+1}{\underline{\lambda}}\bigg\}. \]
Then $\left(A_k,B_k\right)\in\partial H_p(z_k)$ and $\norm{\left(A_k,B_k\right)}\leq M_{2,k}\norm{z_k-z_{k-1}}\leq M_2\norm{z_k-z_{k-1}}$.
\end{lemma}
\proof Invoking~(\ref{Algorithm: x_{k+1}}), one has $$0\in \partial f(x_{k+1})+\nabla g(x_k)+\frac{\alpha_k(x_{k-1}-x_k)}{\lambda_k}+\frac{\nabla h(x_{k+1})-\nabla h(x_{k})}{\lambda_k},$$ which implies that $u_k\in\partial f(x_{k+1})$. It then follows from Lemma~\ref{lem: subdifferential H} that $\left(A_k,B_k\right)\in\partial H_p(z_k)$.  Moreover,
\begin{align*}
\norm{u_k+\nabla g(x_{k+1})}&\leq\norm{\frac{1}{\lambda_k}\left(\nabla h(y_k)-\nabla h(x_{k+1})\right)+\nabla g(x_{k+1})-\nabla g(x_k)}\\
&+2p\norm{x_{k+1}-x_k}+\frac{\alpha_k}{\lambda_k}\norm{x_k-x_{k-1}}\\
&\leq\frac{L_{\nabla h}}{\lambda_k}\norm{y_k-x_k+x_k-x_{k+1}}+\left(L_{\nabla g}+2p \right)\norm{x_{k+1}-x_k}+\frac{1}{\lambda_k}\norm{x_k-x_{k-1}}\\
&\leq \frac{L_{\nabla h}L_{\nabla g}\lambda_{k-1}+1}{\lambda_k}\norm{x_k-x_{k-1}}+\left(\frac{L_{\nabla h}}{\lambda_k}+L_{\nabla g}+2p\right)\norm{x_{k+1}-x_k}.
\end{align*}
Then
\begin{align*}
\norm{\left(A_k,B_k\right)}&\leq\norm{A_k}+\norm{B_k}\leq\norm{u_k+\nabla g(x_{k+1})}+4p\norm{(x_{k+1}-x_k)}\\
&\leq\frac{L_{\nabla h}L_{\nabla g}\lambda_{k-1}+1}{\lambda_k}\norm{x_k-x_{k-1}}+\left(\frac{L_{\nabla h}}{\lambda_k}+L_{\nabla g}+6p\right)\norm{x_{k+1}-x_k},\\
&\leq M_{2,k}\norm{z_k-z_{k-1}}\leq M_{2}\norm{z_k-z_{k-1}},
\end{align*}
as desired.\endproof

Having established basic properties of $\big(H_{p_k}(x_{k+1},x_k) \big)$, we now show its convergence. Denote by $\omega(z_{-1})$ the set of all limit points of $(z_k)_{k\in\N^*}$.

\begin{theorem}\label{thm: function value convergence} Let $(x_k)_{k\in\N^*}$ be a sequence generated by BiFRB and define $z_k=(x_{k+1},x_k)$ for $k\in\N^*$. Suppose that Assumption~\ref{assmption: p_k and M_1,k} holds and $(z_k)_{k\in\N^*}$ is bounded. Let $(z_{k_l})_{l\in\N}$ be a subsequence such that $z_{k_l}\to z^*$ for some $z^*=(x^*,y^*)$ as $l\to\infty$. Then the following hold:
	
(i) We have $\lim_{l\to\infty}H_{p_{k_l}}(z_{k_l})=F(x^*)$. In fact, $\lim_{k\to\infty}H_{p_{k}}(z_{k})=F(x^*)$. Consequently $\lim_{k\to\infty}H_{\bar{p}}(z_{k})=F(x^*)$, where $\bar{p}$ is given in Assumption~\ref{assmption: p_k and M_1,k}(ii).

(ii) The limit point $z^*=(x^*,y^*)$ satisfies $x^*=y^*$ and $0\in\partial H_{\bar{p}}(x^*,y^*)$, which implies that $0\in\partial F(x^*)$.

(iii) The set $\omega(z_{-1})$ is nonempty, connected, and compact, on which the function $H_{\bar{p}}$ is constant $F(x^*)$. Moreover, $\lim_{k\to\infty}\dist(z_k,\omega(z_{-1}))=0$.
\end{theorem}
\proof (i) For every $k\in\N$, appealing to the iterative scheme~(\ref{Algorithm: x_{k+1}}), we have
\begin{align}\label{fom: x}
f(x_{k+1})&\leq f(x^*)+\ip{x^*-x_{k+1}}{\nabla g(x_k)}+\frac{\alpha_k}{\lambda_k}\ip{x^*-x_{k+1}}{x_{k-1}-x_k}\nonumber \\
&+\frac{1}{\lambda_k}D_h(x^*,y_k)-\frac{1}{\lambda_k}D_h(x_{k+1},y_k)\nonumber\\
&\leq f(x^*)+\ip{x^*-x_{k+1}}{\nabla g(x_k)}+\frac{\alpha_k}{\lambda_k}\ip{x^*-x_{k+1}}{x_{k-1}-x_k}\nonumber\\
&+\frac{L_{\nabla h}-\sigma}{2\lambda_k}\norm{x^*-y_k}^2-\frac{\sigma}{2\lambda}\big(\norm{x_{k+1}-y_k}^2-\norm{x^*-y_k}^2\big)\nonumber\\
&=f(x^*)+\ip{x^*-x_{k+1}}{\nabla g(x_k)}+\frac{\alpha_k}{\lambda_k}\ip{x^*-x_{k+1}}{x_{k-1}-x_k}\nonumber\\
&+\frac{L_{\nabla h}-\sigma}{2\lambda_k}\norm{x^*-y_k}^2-\frac{\sigma}{2\lambda_k}\big(\norm{x_{k+1}-x^*}+2\ip{x_k-x^*}{x^*-y_k}\big).
\end{align}
Recall that $x_{k_l}\to x^*$ as $l\to\infty$. Then Lemma~\ref{lem: merit function descent}(ii) and the triangle inequality imply that $x_{k_l-1}\to x^*$ and
\[\norm{x^*-y_{k_l-1}}\leq\norm{x^*-x_{k_l-1}}+\lambda_{k_l-2}L_{\nabla g}\norm{x_{k_l-2}-x_{k_l-1}}\to0,\text{ as }l\to\infty,\]
where the first inequality holds due to~(\ref{Algorithm: y_k}).
Taking the above limits into account, replacing the index $k+1$ by $k_l$ in~(\ref{fom: x}) and passing to the limit, one gets
\[\limsup_{l\to\infty}f(x_{k_l})\leq f(x^*)\leq \liminf_{l\to\infty} f(x_{k_l}),\]
where the last inequality holds thanks to the lower semicontinuity of $f$, which implies that~$\lim_{l\to\infty} F(x_{k_l})=F(x^*)$. Note that the sequence $(H_{p_k}(z_k))_{k\in\N}$ is convergent. Then one concludes that
$$\lim_{k\to\infty} H_{p_k}(z_k)=\lim_{l\to\infty}H_{p_{k_l}}(z_{k_l})=\lim_{l\to\infty}\big(F(x_{k_l})+p_{k_l}\norm{x_{k_l}-x_{k_l-1}}^2 \big)=F(x^*),$$
where $\lim_{l\to\infty}p_{k_l}\norm{x_{k_l}-x_{k_l-1}}^2=0$ due to Assumption~\ref{assmption: p_k and M_1,k}(ii) and Lemma~\ref{lem: merit function descent}(ii).  Finally, it is easy to see  that $0\leq(\bar{p}-p_k)\norm{x_{k}-x_{k-1}}^2\leq \bar{p}\norm{x_{k}-x_{k-1}}^2\to0$ as $k\to\infty$ under Assumption~\ref{assmption: p_k and M_1,k}(ii). Thus $\lim_{k\to\infty}H_{\bar{p}}(z_k)=\lim_{k\to\infty}\big(H_{p_k}(z_k)+(\bar{p}-p_k)\norm{x_{k}-x_{k-1}}^2 \big)=F(x^*)$.

(ii) By Lemma~\ref{lem: merit function descent}(ii) and the triangle inequality, subsequences $(x_{k_l})_{l\in\N}$ and $(x_{k_{l+1}})_{l\in\N}$ have the same limit. Therefore $x^*=y^*$. Combining Lemma~\ref{lem: subdifferential bound} with $p=\bar{p}$, Lemma~\ref{lem: merit function descent}(ii), statement(i), and the outer semicontinuity ot $z\mapsto \partial H_{\bar{p}}(z)$, we conclude that $0\in\partial H_{\bar{p}}(z^*)$, which further implies $0\in\partial F(x^*)$.

(iii) Apply a similar proof as in~\cite[Theorem 3.6(iii)]{wang2021malitsky}.\endproof

The result below provides a sufficient condition to the critical boundedness assumption in Theorem~\ref{thm: function value convergence}.

\begin{theorem}\label{thm: bounded sequence} Let $(x_k)_{k\in\N^*}$ be a sequence generated by BiFRB and define $z_k=(x_{k+1},x_k)$ for $k\in\N^*$. Suppose that Assumption~\ref{assmption: p_k and M_1,k}(i) holds. If $f+g$ is coercive (or level-bounded), then the sequence $(x_k)_{k\in\N^*}$ is bounded, so is $(z_k)_{k\in\N^*}$.
\end{theorem}
\proof Assume without loss of generality that $M_{1,k}>0$ and $p_k\geq0$ for $k\in\N$.  Then, appealing to Lemma~\ref{lem: merit function descent}, one has $H_{p_k}(z_k)\leq H_{p_{-1}}(z_{-1})$ for $k\in\N$. Therefore
\begin{equation*}
 F(x_{k+1})=H_{p_k}(z_k)-p_k\norm{x_{k+1}-x_k}^2\leq H_{p_k}(z_k)\leq H_{p_{-1}}(z_{-1}).
\end{equation*}
If $(x_k)_{k\in\N}$ was unbounded, then the above inequality would yield a contradiction due to the coercivity of $f+g$ (or level-boundedness). \endproof
\subsection{Abstract sequential convergence}\label{sec:Sequential convergence}

In this section, we establish global sequential convergence of BiFRB to a stationary point by using the generalized concave KL property (recall Definition~\ref{def: g-concave KL}). For $\varepsilon>0$ and nonempty set $\Omega\subseteq\Rn$, we define $\Omega_\varepsilon=\{x\in\Rn:\dist(x,\Omega)<\varepsilon\}$.

\begin{lemma}\label{lemma: Uniformize the g-concave KL}\emph{\cite[Lemma 4.4]{wang2020}} Let $f:\Rn\rightarrow\overline{\R}$ be proper lsc and let $\mu\in\R$. Let $\Omega\subseteq\dom\partial f$ be a nonempty compact set on which $f(x)=\mu$ for all $x\in\Omega$. The following statements hold:
	
	(i) Suppose that $f$ satisfies the pointwise generalized concave KL property at each $x\in\Omega$. Then there exist $\varepsilon>0,\eta\in(0,\infty]$ and $\varphi\in\Phi_\eta$ such that $f$ has the setwise generalized concave KL property on $\Omega$ with respect to $U=\Omega_\varepsilon$, $\eta$ and $\varphi$.
	
	(ii) Set $U=\Omega_\varepsilon$ and define $h:(0,\eta)\rightarrow\R_+$ by \[h(s)=\sup\big\{\dist^{-1}\big(0,\partial f(x)\big):x\in U\cap[0<f-\mu<\eta],s\leq f(x)-\mu\big\}.\]
	Then the function $\tilde{\varphi}:[0,\eta)\rightarrow\R_+$,
	\[t\mapsto\int_0^th(s)ds,~\forall t\in(0,\eta),\]
	and $\tilde{\varphi}(0)=0$, is well-defined and belongs to $\Phi_\eta$. The function $f$ has the setwise generalized concave KL property on $\Omega$ with respect to $U$, $\eta$ and $\tilde{\varphi}$. Moreover,
	\[\tilde{\varphi}=\inf\big\{\varphi\in\Phi_\eta:\text{$\varphi$ is a concave desingularizing function of $f$ on $\Omega$ with respect to $U$ and $\eta$}\big\}.\]
	We say $\tilde{\varphi}$ is the exact modulus of the setwise generalized concave KL property of $f$ on $\Omega$ with respect to $U$ and $\eta$.
\end{lemma}

\begin{theorem}[global convergence of BiFRB]\label{thm: finite length}  Let $(x_k)_{k\in\N^*}$ be a sequence generated by BiFRB and define $z_k=(x_{k+1},x_k)$ for $k\in\N^*$. Suppose that Assumption~\ref{assmption: p_k and M_1,k} holds and $(z_k)_{k\in\N^*}$ is bounded. Suppose further that the dominating BiFRB merit function $H_{\bar{p}}$ (recall Definition~\ref{def: merit}) has the generalized concave KL property on $\omega(z_{-1})$. Then the following hold:
	
	(i) The sequence $(z_k)_{k\in\N^*}$ is Cauchy and has finite length. To be specific, there exist index $k_0\in\N$, $\varepsilon>0$ and $\eta\in(0,\infty]$ such that for $i\geq k_0+1$,
	\begin{equation}\label{fom: finite lenght upper bound}
		\sum_{k= i }^ \infty\norm{z_{k+1}-z_k}\leq\norm{z_i-z_{i-1}}+C\tilde{\varphi}\left(H_{p_i}(z_{ i })-F(x^*)\right),
	\end{equation}
	where
	$$C=\frac{2\sqrt{2}\max\bigg\{\frac{L_{\nabla h}}{\underline{\lambda}}+L_{\nabla g}+6\bar{p},  \frac{L_{\nabla h}L_{\nabla g}\overline{\lambda}+1}{\underline{\lambda}}\bigg\}}{\liminf_{k\to\infty}M_{1,k}},$$ 
	and $\tilde{\varphi}$ is the exact modulus associated with the setwise generalized concave KL property of $H_{\bar{p}}$ on $\omega(z_{-1})$ with respect to $\varepsilon$ and $\eta$.
	
	(ii) The sequence $(x_k)_{k\in\N^*}$ has finite length and converges to some $x^*$ with $0\in\partial F(x^*)$.
\end{theorem}
\proof By the boundedness assumption, assume without loss of generality that $z_k\to z^*=(x^*,y^*)\in\Rn\times\Rn$. Then Theorem~\ref{thm: function value convergence} implies that $x^*=y^*$ and $H_{p_k}(z_k)\to F(x^*)$ as $k\to\infty$. By Lemma~\ref{lem: merit function descent}, assume without loss of generality that $\big(H_{p_k}(z_k) \big)_{k\in\N}$ is decreasing. Consider two cases:

Case 1: Suppose that there exists $k_0$ such that $F(z^*)=H_{p_{k_0}}(z_{k_0})$. Then Lemma~\ref{lem: merit function descent}(i) implies that $z_{k_0+1}=z_{k_0}$ and $x_{k_0+1}=x_{k_0}$. The desired results then follows from a simple induction.

Case 2: Assume that $F(z^*)<H_{p_k}(z_k)$ for all $k$.  By Theorem~\ref{thm: function value convergence}(iii), the function $H_{\bar{p}}$ is constant $F(x^*)$ on $\omega(z_{-1})$. Then Lemma~\ref{lemma: Uniformize the g-concave KL} yields $\varepsilon>0$ and $\eta>0$ such that $H_{\bar{p}}$ has the setwise generalized concave KL property on $\omega(z_{-1})$ with respect to $\varepsilon>0$ and $\eta>0$ and the associated exact modulus $\tilde{\varphi}$. Appealing to Theorem~\ref{thm: function value convergence} again, there exists $k_1$ such that $z_k\in[0<H_{\bar{p}}-H_{\bar{p}}(z^*)<\eta]$ for $k>k_1$ and $k_2$ such that $\dist(z_k,\omega(z_{-1}))<\varepsilon$ for all $k>k_2$. Put $k_0=\max(k_1,k_2)$. Then for $k>k_0$,
\begin{align*}
	(\tilde{\varphi})_-^\prime\left(H_{p_k}(z_k)-F(x^*) \right)\cdot\dist(0,\partial H_{\bar{p}}(z_k))\geq (\tilde{\varphi})_-^\prime\left(H_{\bar{p}}(z_k)-F(x^*) \right)\cdot\dist(0,\partial H_{\bar{p}}(z_k))\geq1,
\end{align*}
where the first inequality holds because $H_{p_k}(z_k)\leq H_{\bar{p}}(z_k)$ and the last one follows from the generalized concave KL inequality (recall Definition~\ref{def: g-concave KL}). Invoking Lemma~\ref{lem: subdifferential bound} yields
\begin{equation}\label{d}
	M_{2,k} (\tilde{\varphi})_-^\prime\left(H_{p_k}(z_k)-F(x^*) \right)\cdot\norm{z_k-z_{k-1}}\geq1,
\end{equation}
where $M_{2,k}=M_{2,k}(\bar{p})$. For simplicity, define  $$\Delta_{k,k+1}=\tilde{\varphi}\left(H_{p_k}(z_k)-H(z^*)\right)-\tilde{\varphi}\left(H_{p_{k+1}}(z_{k+1})-H(z^*)\right).$$ By the concavity of $\tilde{\varphi}$ and~(\ref{d}), one has
\begin{align}\label{cc}
	\Delta_{k,k+1}&=\tilde{\varphi}\left(H_{p_k}(z_k)-H(z^*)\right)-\tilde{\varphi}\left(H_{p_{k+1}}(z_{k+1})-H(z^*)\right)\nonumber \\
	&\geq(\tilde{\varphi})_-^\prime\left(H_{p_k}(z_k)-H(z^*)\right)\cdot[H_{p_k}(z_k)-H_{p_{k+1}}(z_{k+1})]\nonumber\\
	&\geq\frac{H_{p_k}(z_k)-H_{p_{k+1}}(z_{k+1})}{M_{2,k}\norm{z_k-z_{k-1}}},
\end{align}
Applying Lemma~\ref{lem: merit function descent} to~(\ref{cc}) yields
\begin{equation*}
	\Delta_{k,k+1}\geq\frac{\norm{z_{k+1}-z_k}^2}{C_k\norm{z_k-z_{k-1}}},
\end{equation*}
where $C_k=M_{2,k}/M_{1,k}$. On one hand, by Assumption~\ref{assmption: p_k and M_1,k}(i), assume without loss of generality that $M_{1,k}\geq \liminf_{k\to\infty}M_{1,k}/2$ for all $k$. On the other hand, $M_{2,k}$ satisfies
$$M_{2,k}\leq \sqrt{2}\max\bigg\{\frac{L_{\nabla h}}{\underline{\lambda}}+L_{\nabla g}+6\bar{p},  \frac{L_{\nabla h}L_{\nabla g}\overline{\lambda}+1}{\underline{\lambda}}\bigg\}.$$
Then $C_k\leq C$ and
\begin{align}\label{xx}
	2\norm{z_{k+1}-z_k}\leq 2\sqrt{C_k\Delta_{k,k+1}\norm{z_k-z_{k-1}}}\leq C\Delta_{k,k+1}+\norm{z_k-z_{k-1}}.
\end{align}
Pick $i\geq k_0+1$. Summing~(\ref{xx}) from $i$ to an arbitrary $j>i$ yields
\begin{align*}
	&~~~~2\sum_{k= i }^ j\norm{z_{k+1}-z_k}\leq\sum_{k= i }^ j\norm{z_k-z_{k-1}}+C\sum_{k= i }^ j \Delta_{k,k+1}\\
	&\leq\sum_{k=i}^ j\norm{z_{k+1}-z_k}+\norm{z_{ i }-z_{i-1}}+C\tilde{\varphi}\left( H_{p_i}(z_{ i })-F(x^*)\right)-C\tilde{\varphi}\left( H_{p_{j+1}}(z_{ j+1})-H(z^*)\right)\\
	&\leq \sum_{k= i }^ j\norm{z_{k+1}-z_k}+\norm{z_{ i }-z_{i-1}}+C\tilde{\varphi}\left(  H_{p_i}(z_{ i })-F(x^*)\right),
\end{align*}
implying that
\begin{equation*}
	\sum_{k= i }^ j\norm{z_{k+1}-z_k}\leq\norm{z_i-z_{i-1}}+C\tilde{\varphi}\left(H_{p_i}(z_{ i })-F(x^*)\right),
\end{equation*}
from which~(\ref{fom: finite lenght upper bound}) readily follows. Finally, notice that
\begin{equation}\label{formula: towards Cauchy}
	\norm{z_{j}-z_i}\leq\sum_{k= i }^ j\norm{z_{k+1}-z_k}\leq\norm{z_i-z_{i-1}}+C\tilde{\varphi}\left(H_{p_i}(z_{ i })-F(x^*)\right).
\end{equation}
Recall from Theorem~\ref{thm: function value convergence}(i) and Lemma~\ref{lem: merit function descent}(ii) that $H_{p_i}(z_i)\to F(x^*)$ and $\norm{z_i-z_{i-1}}\to0$ as $i\to\infty$. Thus, passing the index $i$ to infinity, we conclude from inequality~(\ref{formula: towards Cauchy}) that $(z_k)_{k\in\N}$ is Cauchy. Invoking Theorem~\ref{thm: function value convergence}(ii) completes the proof.

(ii) The statement follows from the definition of $(z_k)_{k\in\N^*}$ and Theorem~\ref{thm: function value convergence}(ii).\endproof

The generalized concave KL assumption in Theorem~\ref{thm: finite length} can be easily satisfied by semialgebraic functions.

\begin{corollary}\label{cor: convergence on semialgebraic functions} Let $(x_k)_{k\in\N^*}$ be a sequence generated BiFRB. Assume that $(x_k)_{k\in\N^*}$ is bounded and Assumption~\ref{assmption: p_k and M_1,k} holds. Suppose further that both $f$ and $g$ are semialgebraic. Then $(x_k)_{k\in\N^*}$ converges to some $x^*$ with $0\in\partial F(x^*)$ and has finite length property.
\end{corollary}
\proof Let $\bar{p}$ be as in Assumption~\ref{assmption: p_k and M_1,k}. Then $(x,y)\mapsto\bar{p}\norm{x-y}^2$ is semialgebraic. In turn, the dominating BiFRB function $H_{\bar{p}}$ (recall Definition~\ref{def: merit}) is semialgebraic by the semialgebraic functions sum rule; see, e.g.,~\cite[Section 4.3]{Attouch2010}. The desired result then follows immediately from Theorem~\ref{thm: finite length}.~\endproof

\section{Stepsize strategies}\label{sec:Stepsize strategies }
Following abstract convergence properties of BiFRB, we now exploit conditions such that the important Assumption~\ref{assmption: p_k and M_1,k} holds. 

To furnish Assumption~\ref{assmption: p_k and M_1,k}, it suffices to establish the existence of suitable $p_{-1}$ under desirable stepsize rules. Unlike most literature that emphasizes on explicit merit function parameters~\cite{boct2016inertial,boct2016inertiala,wang2021malitsky,ochs2014ipiano},
our $(p_k)_{k\in\N^*}$ and $(M_{1,k})_{k\in\N}$ are implicit. This turns out to be instrumental in designing less restrictive stepsize rules.
  In the remainder of this section, unless specified otherwise,
\begin{equation}\label{def: abc}
	a=(L_{\nabla h}-\sigma)L_{\nabla g}^2, b=\sigma, \text{ and }c=L_{\nabla g}.
\end{equation}
Consequently, the sequence $(p_k)_{k\in\N^*}$ and $(M_{1,k})_{k\in\N^*}$ given by~(\ref{def: p_k}) and~(\ref{def:M_1,k}) satisfy
\begin{align}\label{simple pk}
	p_{k}=\frac{a}{2}\left(\frac{\lambda_{k-1}^2}{\lambda_k}\right)+\frac{1}{2}\left(\frac{b}{\lambda_k}-c\right)-p_{k-1},~M_{1,k}=p_{k-1}-\frac{\alpha_k+b c\lambda_{k-1}}{2\lambda_k}-\frac{a\lambda_{k-1}^2}{2\lambda_k}.
\end{align}

The lemma below asserts that $(p_k)_{k\in\N^*}$ and $(M_{1,k})_{k\in\N}$ are governed by stepsizes and an initial input $p_{-1}$, which will be useful soon.
\begin{lemma}\label{lem: p_k} Let $(\lambda_k)_{k\in\N^*}$ be a sequence of nonzero real numbers and let $a,b,c\in\R$. Set $p_{-1}\in\R$ and define $(p_k)_{k\in\N^*}$ as in~(\ref{simple pk}).  Then for $k\in\N$,
\begin{align*}
&p_{2k}=\frac{a}{2}\sum_{i=1}^k \left(\frac{\lambda_{2i-1}^2}{\lambda_{2i}}-\frac{\lambda_{2i-2}^2}{\lambda_{2i-1}}\right)+\frac{b}{2}\sum_{i=1}^k\left(\frac{1}{\lambda_{2i}}-\frac{1}{\lambda_{2i-1}}\right)+\frac{a\lambda^2_{-1}+b}{2\lambda_0}-\frac{c}{2}-p_{-1},\\
&p_{2k+1}=\frac{a}{2}\sum_{i=0}^k\left(\frac{\lambda_{2i}^2 }{\lambda_{2i+1}}-\frac{\lambda_{2i-2}^2}{\lambda_{2i-1}}\right)+\frac{b}{2}\sum_{i=0}^k\left(\frac{1}{\lambda_{2i+1}}-\frac{1}{\lambda_{2i}}\right)+p_{-1},
\end{align*}
where we adopt the convention $\sum_{i=1}^0a_i=0$ for any sequence $(a_k)_{k\in\N}$.
\end{lemma}
\proof Let $k\in\N\backslash\{0\}$. Then
\begin{align*}
p_{2k}&=\frac{a}{2}\left(\frac{\lambda_{2k-1}^2}{\lambda_{2k}}\right)+\frac{1}{2}\left(\frac{b}{\lambda_{2k}}-c\right)-p_{2k-1}=\frac{a}{2}\left(\frac{\lambda_{2k-1}^2}{\lambda_{2k}}-\frac{\lambda_{2k-2}^2}{\lambda_{2k-1}}\right)+\frac{b}{2}\left(\frac{1}{\lambda_{2k}}-\frac{1}{\lambda_{2k-1}}\right)+p_{2k-2}\\
&=\frac{a}{2}\sum_{i=0}^k \left(\frac{\lambda_{2i-1}^2}{\lambda_{2i}}-\frac{\lambda_{2i-2}^2}{\lambda_{2i-1}}\right)+\frac{b}{2}\sum_{i=1}^k\left(\frac{1}{\lambda_{2i}}-\frac{1}{\lambda_{2i-1}}\right)+p_0\\
&=\frac{a}{2}\sum_{i=1}^k \left(\frac{\lambda_{2i-1}^2}{\lambda_{2i}}-\frac{\lambda_{2i-2}^2}{\lambda_{2i-1}}\right)+\frac{b}{2}\sum_{i=1}^k\left(\frac{1}{\lambda_{2i}}-\frac{1}{\lambda_{2i-1}}\right)+\frac{a\lambda^2_{-1}+b}{2\lambda_0}-\frac{c}{2}-p_{-1},
\end{align*}
and
\begin{align*}
p_{2k+1}&=\frac{a}{2}\left(\frac{\lambda_{2k}^2}{\lambda_{2k+1}} \right)+\frac{1}{2}\left(\frac{b}{\lambda_{2k+1}}-c\right)-p_{2k}\\
&=\frac{a}{2}\left(\frac{\lambda_{2k}^2}{\lambda_{2k+1}} \right)+\frac{1}{2}\left(\frac{b}{\lambda_{2k+1}}-c\right)-\frac{a}{2}\sum_{i=1}^k \left(\frac{\lambda_{2i-1}^2}{\lambda_{2i}}-\frac{\lambda_{2i-2}^2}{\lambda_{2i-1}}\right)-\frac{b}{2}\sum_{i=1}^k\left(\frac{1}{\lambda_{2i}}-\frac{1}{\lambda_{2i-1}}\right)\\
&-\frac{a\lambda^2_{-1}+b}{2\lambda_0}+\frac{c}{2}+p_{-1}\\
&=\frac{a}{2}\sum_{i=0}^k\left(\frac{\lambda_{2i}^2 }{\lambda_{2i+1}}-\frac{\lambda_{2i-2}^2}{\lambda_{2i-1}}\right)+\frac{b}{2}\sum_{i=0}^k\left(\frac{1}{\lambda_{2i+1}}-\frac{1}{\lambda_{2i}}\right)+p_{-1},
\end{align*}
from which the desired identities hold. \endproof

\subsection{Non-Euclidean case}\label{sec: non-Euclidean stepsize}
Now we provide a fixed stepsize rule for Algorithm~\ref{BiFRB}. In contrast to similar stepsize results in the literature, we shall see that Assumption~\ref{assmption: p_k and M_1,k} furnishes explicit stepsize bounds that is independent of inertial parameters. In turn, we answer a question of Malitsky and Tam \cite[Section 7]{malitsky2020FRB};
see Remark~\ref{rem: parameter condition advantage}.


\begin{proposition}[fixed stepsize]\label{prop: non-Euclidean fixed stepsize} Let $(\lambda_k)_{k\in\N^*}$ be a constant BiFRB sequence with $\lambda_k=\lambda>0$ for all $k$. Suppose that $\sigma>2$, $(L_{\nabla h}-\sigma)\sigma>1/4$ and $0\leq\alpha_k<\min(1,\sigma/2)=1$. Define
	\[ \lambda^*=\frac{\sqrt{(2bc+c)^2+4a(b-2)}-2bc-c}{2a}>0, \]
	where $a,b,c$ are specified by~(\ref{def: abc}). If $\lambda<\min\big\{ \lambda^*,(\sigma-1)/[(\sigma+1)L_{\nabla g}]\big\}$, then there exists $p_{-1}$ such that $(p_k)_{k\in\N^*}$ and $(M_{1,k})_{k\in\N}$ generated by~(\ref{def: p_k}) and~(\ref{def:M_1,k}) satisfy Assumption~\ref{assmption: p_k and M_1,k}.
\end{proposition}
\proof Applying Lemma~\ref{lem: p_k} to~(\ref{simple pk}) yields that $p_k=p_{-1}$ if $k$ is odd; $p_k= a\lambda/2 +\left(b/\lambda-c\right)/2-p_{-1}$ if $k$ is even; and
\begin{align*}
	M_{1,k}=\begin{cases*}
		p_{-1}-\frac{\alpha_{k}}{2\lambda}-\frac{bc}{2}-\frac{a\lambda}{2},&\text{if $k$ is even,}\\
		\frac{b-\alpha_k}{2\lambda}-\frac{(b+1)c}{2}-p_{-1},&\text{if $k$ is odd.}
	\end{cases*}
\end{align*}
Then\begin{align}\label{formula: non-Euclidean pk}
\min\bigg\{p_{-1},\frac{a}{2}\lambda+\frac{1}{2}\left(\frac{b}{\lambda}-c \right)-p_{-1}\bigg\}\leq p_k\leq\max\bigg\{p_{-1},\frac{a}{2}\lambda+\frac{1}{2}\left(\frac{b}{\lambda}-c \right)-p_{-1}\bigg\},
\end{align}
\begin{equation}\label{formula: non-Euclidean M1,k}
	\liminf_{k\to\infty}M_{1,k}\geq \min\bigg\{p_{-1}-\frac{1}{2\lambda}-\frac{bc}{2}-\frac{a\lambda}{2},\frac{b-1}{2\lambda}-\frac{(b+1)c}{2}-p_{-1}\bigg\}.
\end{equation}

To furnish Assumption~\ref{assmption: p_k and M_1,k}, we will prove that
\begin{align}
	&\frac{a}{2}\lambda+\frac{1}{2}\left(\frac{b}{\lambda}-c \right)>0,\label{fom: fixed stepsize 1}\\
&\frac{1}{2\lambda}+\frac{bc}{2}+\frac{a\lambda}{2}<\frac{b-1}{2\lambda}-\frac{(b+1)c}{2},\label{fom: fixed stepsize 2}\\
&\frac{1}{2\lambda}+\frac{bc}{2}+\frac{a\lambda}{2}<\frac{a\lambda}{2}+\frac{b}{2\lambda}-\frac{c}{2}.\label{fom: fixed stepsize 3}
\end{align}
Combining inequalities~(\ref{fom: fixed stepsize 1})--(\ref{fom: fixed stepsize 3}) yields the existence of some $p_{-1}$ such that
\[0<p_{-1}<\frac{a}{2}\lambda+\frac{1}{2}\left(\frac{b}{\lambda}-c \right)\text{ and }\frac{1}{2\lambda}+\frac{bc}{2}+\frac{a\lambda}{2}<p_{-1}<\frac{b-1}{2\lambda}-\frac{(b+1)c}{2},\]
which together with~(\ref{formula: non-Euclidean pk}) and~(\ref{formula: non-Euclidean M1,k}) entails Assumption~\ref{assmption: p_k and M_1,k}. First, we establish~(\ref{fom: fixed stepsize 1}). Observe that
\[\frac{a}{2}\lambda+\frac{1}{2}\left(\frac{b}{\lambda}-c \right)>0\Leftrightarrow a\lambda^2-c\lambda+b>0.\]
By assumption, $a>0$ and $c^2-4ab=L^2_{\nabla g}\big(1-4(L_{\nabla h}-\sigma)\sigma \big)<0$. Then the quadratic function $\lambda\mapsto a\lambda^2-c\lambda+b$ has no root and is therefore strictly positive, implying~(\ref{fom: fixed stepsize 1}). Next we justify~(\ref{fom: fixed stepsize 2}), which amounts to $-a\lambda^2-(2bc+c)\lambda+b-2>0$. Notice that the discriminant $(2bc+c)^2+4a(b-2)>0$ and product of roots $(b-2)/(-a)<0$. Then the quadratic function $\lambda\mapsto-a\lambda^2-(2bc+c)\lambda+b-2$ has two real roots with opposite signs. In turn, $$(\forall~0<\lambda< \lambda^*)-a\lambda^2-(2bc+c)\lambda+b-2>0,$$ furnishing inequality~(\ref{fom: fixed stepsize 2}). Finally, the assumption $\lambda<(\sigma-1)/[(\sigma+1)L_{\nabla g}]=(b-1)/[(b+1)c]$ implies~(\ref{fom: fixed stepsize 3}) after rearranging.\endproof
\begin{remark}\label{rem: parameter condition advantage} (i)~(Comparison to known results) Compared to a result of Bo{\c{t}} and Csetnek~\cite[Lemma 3.3]{boct2016inertial} with explicit 
merit function parameters but implicit stepsize bounds, our stepsize is explicit and independent of inertial parameters; see~\cite[Remark 7]{boct2016inertiala} for a similar restrictive result and~\cite[Corollary 4.4]{malitsky2020FRB} for a result with dependent stepsize bound and inertial parameter in the presence of convexity; see also~\cite[Table 1]{ochs2018local}.

(ii)~(Regarding a question of Malitsky and Tam) In~\cite[Section 7]{malitsky2020FRB}, Malitsky and Tam posted a question regarding whether FRB can be adapted to incorporate a Nesterov-type acceleration. With inertial parameter $\alpha_k\in[0,1)$ independent of stepsize $\lambda$ in Proposition~\ref{prop: non-Euclidean fixed stepsize}, as opposed to a fixed inertial parameter $\alpha\in[0,1/3)$ dependent of stepsize $\lambda$ in~\cite[Corollary 4.4]{malitsky2020FRB}, this question is resolved. Indeed, the Nesterov acceleration scheme corresponds to Proposition~\ref{prop: non-Euclidean fixed stepsize} with \[(\forall k\in\N^*)~\alpha_k=\frac{t_k-1}{t_{k+1}},\]
where $(\forall k\in\N^*)~t_{k+1}=(1+\sqrt{1+4t_k^2})/2$ and $t_{-1}=1$; see, e.g.,~\cite{nesterov1983method}. 

(iii) The assumption that $(L_{\nabla h}-\sigma)\sigma>1/4$ is not demanding. For instance, let $\alpha,\beta>0$ be such that $\alpha\beta>1/4$ and let $u:\Rn\to\R$ be a function with $1$-Lipschitz gradient. Then $h=\alpha u+\beta\norm{\cdot}^2/2$ satisfies $(L_{\nabla h}-\sigma)\sigma>1/4$.
\end{remark}

\subsection{Euclidean case}\label{sec: Euclidean stepsize}
Considering kernel $h=\norm{\cdot}^2/2$, we now turn to stepsize strategies of Algorithm~\ref{iFRB}.

\begin{proposition}[dynamic stepsize] Let $(\lambda_k)_{k\in\N^*}\subseteq[\underline{\lambda},\overline{\lambda}]$ be a iFRB stepsize sequence and let $(p_k)_{k\in\N^*}$ be the sequence generated by~(\ref{simple pk}) with $a=0$, $b=1$ and $0\leq\alpha_k<\bar{\alpha}$ for some $\bar{\alpha}<1/2$. Let $(a_k)_{k\in\N}$ be a sequence of positive real numbers with $\sum_{k=0}^\infty a_k<\infty$ and let $0<\varepsilon<(1-2\bar{\alpha})/(3L_{\nabla g})$. Suppose that the following hold:

(i) The stepsize lower bound $\underline{\lambda}=\varepsilon$ and upper bound $\overline{\lambda}<\varepsilon/(2\bar{\alpha}+3\varepsilon L_{\nabla g} )$.

(ii) For all $k\in\N$, $0\leq(1/\lambda_k)-(1/\lambda_{k-1})\leq a_k$. 	
\vspace{1mm}

Then there exists $p_{-1}$ such that $(p_k)_{k\in\N^*}$ and $(M_{1,k})_{k\in\N}$ generated by~(\ref{def: p_k}) and~(\ref{def:M_1,k}) satisfy Assumption~\ref{assmption: p_k and M_1,k}.
\end{proposition}
\proof First notice that assumption(i) is not contradictory because $\varepsilon<(1-2\bar{\alpha})/(3L_{\nabla g})\Leftrightarrow2\bar{\alpha}+3\varepsilon L_{\nabla g}<1\Leftrightarrow\varepsilon< \varepsilon/(2\bar{\alpha}+3\varepsilon L_{\nabla g} )$. Then, by assumption(i), we have
\begin{align*}
\lambda_{0}<\frac{\varepsilon}{2\bar{\alpha}+3\varepsilon L_{\nabla g} }\Leftrightarrow\frac{L_{\nabla g}}{2}+\frac{\bar{\alpha}}{2\varepsilon}<\frac{1}{2\lambda_{0}}-L_{\nabla g}+\frac{\bar{\alpha}}{2\varepsilon},
\end{align*}
which implies that there exits $p_{-1}$ such that
\begin{align}\label{fom: euclidean p_{-1}}
\frac{L_{\nabla g}}{2}+\frac{\bar{\alpha}}{2\varepsilon}<p_{-1}<\frac{1}{2\lambda_{0}}-L_{\nabla g}+\frac{\bar{\alpha}}{2\varepsilon}.
\end{align}

According to Lemma~\ref{lem: p_k} and assumption(ii)
\begin{align*}
p_{-1}\leq p_{2k+1}&=\frac{1}{2}\sum_{i=0}^k\left(\frac{1}{\lambda_{2i+1}}-\frac{1}{\lambda_{2i}}\right)+p_{-1}\leq\frac{1}{2}\sum_{k=0}^\infty a_k+p_0 ,\\
\frac{1}{2\lambda_0}-\frac{c}{2}-p_{-1}\leq p_{2k}&=\frac{1}{2}\sum_{i=1}^k\left(\frac{1}{\lambda_{2i}}-\frac{1}{\lambda_{2i-1}} \right)+\frac{1}{2\lambda_0}-\frac{c}{2}-p_{-1}\leq\frac{1}{2}\sum_{k=0}^\infty a_k+\frac{1}{2\lambda_0}-\frac{c}{2}-p_{-1},
\end{align*}
which together with~(\ref{fom: euclidean p_{-1}}) shows that $p_{k}$ is positive and bounded above.  By assumption(ii), $\lambda_k\leq\lambda_{k-1}$ for all $k$. Consequently, the sequence $(\lambda_k)$ is convergent and $\lim_{k\to\infty}\lambda_{k-1}/\lambda_k=1$. It follows from the above inequalities and~(\ref{fom: euclidean p_{-1}}) that
\begin{align*}
	M_{1,2k+1}&=p_{2k+1}-\frac{\alpha_{2k+1}+L_{\nabla g}\lambda_{2k}}{2\lambda_{2k+1}}\geq p_{-1}-\frac{L_{\nabla g}}{2}\left(\frac{\lambda_{2k}}{\lambda_{2k+1}}\right)-\frac{\bar{\alpha}}{2\varepsilon}\to p_{-1}-\frac{L_{\nabla g}}{2}-\frac{\bar{\alpha}}{2\varepsilon}>0,\\
M_{1,2k}&=p_{2k}-\frac{\alpha_{2k}+L_{\nabla g}\lambda_{2k-1}}{2\lambda_{2k}}\\
&\geq \frac{1}{2\lambda_0}-\frac{L_{\nabla g}}{2}\left(\frac{\lambda_{2k-1}}{\lambda_{2k}}+1\right)-\frac{\bar{\alpha}}{2\varepsilon}-p_{-1}\to\frac{1}{2\lambda_0}-L_{\nabla g}-\frac{\bar{\alpha}}{2\varepsilon}-p_{-1}>0.
\end{align*}
Hence $\liminf_{k\to\infty}M_{1,k}>0$.\endproof

The following corollary follows immediately, which provides an inertial variant of~\cite[Lemma 3.2]{wang2021malitsky} together with Lemma~\ref{lem: merit function descent}.
\begin{corollary}[fixed stepsize]\label{cor: Euclidean fixed stepsize} Let $(\lambda_k)_{k\in\N^*}$ be a constant iFRB stepsize sequence with $\lambda_k=\lambda>0$ for all $k$ and let $(p_k)_{k\in\N^*}$ be the sequence generated by~(\ref{simple pk}) with $a=0$, $b=1$ and $0\leq\alpha_k<\bar{\alpha}$ for some $\bar{\alpha}<1/2$. If $\lambda<(1-2\bar{\alpha})/(3L_{\nabla g})$, then there exists $p_{-1}$ such that $(p_k)_{k\in\N^*}$ and $(M_{1,k})_{k\in\N}$ generated by~(\ref{def: p_k}) and~(\ref{def:M_1,k}) satisfy Assumption~\ref{assmption: p_k and M_1,k}.
\end{corollary}
\begin{remark} One recovers the stepsize requirement in~\cite[Lemma 3.2]{wang2021malitsky} by setting $\bar{\alpha}=0$.
\end{remark}

\section{Convergence rates}\label{sec:rate}
After global convergence results, we now turn to convergence rates of both the function value and actual sequence. To this end,
a lemma helps.

\begin{lemma}\emph{\cite[Lemma 10]{boct2020proximal}}\label{lem: general convergence rate} Let $(e_k)_{k\in\N}\subseteq\R_+$ be a decreasing sequence converging $0$. Assume further that there exist natural numbers $k_0\geq l_0\geq1$ such that for every $k\geq k_0$,
	\begin{equation}
		e_{k-l_0}-e_k\geq C_ee_k^{2\theta},
	\end{equation}
where $C_e>0$ is some constant and $\theta\in[0,1)$. Then the following hold:

(i) If $\theta=0$, then $(e_k)_{k\in\N}$ converges in finite steps.

(ii) If $\theta\in\left(0,\frac{1}{2}\right]$, then there exists $C_{e,0}>0$ such that for every $k\geq k_0$,
\[0\leq e_k\leq C_{e,0}Q^k.\]

(iii) If $\theta\in\left(\frac{1}{2},1\right)$, then there exists $C_{e,1}>0$ such that for every $k\geq k_0+l_0$,
\[0\leq e_k\leq C_{e,1}(k-l_0+1)^{-\frac{1}{2\theta-1}}.\]
\end{lemma}

Recall that we say a generalized concave KL function has KL exponent $\theta\in[0,1)$, if an associated desingularizing function $\varphi(t)=c\cdot t^{1-\theta}$ for some $c>0$. We now provide function value convergence rate of BiFRB under KL exponent assumption.

\begin{theorem}[function value convergence rate] \label{thm: function value rate}Let $(x_k)_{k\in\N^*}$ be a sequence generated by BiFRB and define $z_k=(x_{k+1},x_k)$ for $k\in\N^*$. Suppose that all assumptions in Theorem~\ref{thm: finite length} are satisfied and let $x^*$ be the limit given in Theorem~\ref{thm: finite length}(ii). Suppose further that the dominating BiFRB merit function $H_{\bar{p}}$ (recall Definition~\ref{def: merit}) has the generalized concave KL property at $z^*=(x^*,x^*)$ with KL exponent $\theta\in[0,1)$.  Define $e_k=H_{p_k}(z_k)-F(x^*)$ for $k\in\N$. Then $(e_k)_{k\in\N}$ converges to $0$ and the following hold:
	
	(i) If $\theta=0$, then $(e_k)_{k\in\N}$ converges in finite steps.

(ii) If $\theta\in\left(0,1/2\right]$, then there exist $\hat{c}_1>0$ and $\hat{Q}_1\in[0,1)$ such that for $k$ sufficiently large,
\[e_k\leq\hat{c}_1\hat{Q}_1^k.\]

(iii) If $\theta\in\left(1/2,1 \right)$, then there exists $\hat{c}_2>0$ such that for $k$ sufficiently large, \[e_k\leq \hat{c}_2k^{-\frac{1}{2\theta-1}}.\]
\end{theorem}
\proof Recall from Lemma~\ref{lem: merit function descent}(ii) and Theorem~\ref{thm: function value convergence}(i) that the sequence $(e_k)_{k\in\N}$ is decreasing and converges to $0$. By the KL exponent assumption, there exist $c>0$ and $k_0$ such that for $k\geq k_0$
\[\dist\big(0,\partial H_{\bar{p}}(z_k)\big)\geq\frac{\left(H_{p_k}(z_k)-H_{\bar{p}}(z^*) \right)^\theta}{c(1-\theta)}=\frac{e_k^\theta}{c(1-\theta)}.\]
Assume without loss of generality that $M_{1,k}\geq \liminf_{k\to\infty}M_{1,k}/2$. So Lemmas~\ref{lem: merit function descent} and~\ref{lem: subdifferential bound} imply
\begin{align*}
	e_{k-1}-e_k&=H_{p_k-1}(z_{k-1})-H_{p_k}(z_k)\geq M_{1,k}\norm{z_k-z_{k-1}}^2\geq \frac{M_{1,k}}{M_{2,k}^2}\dist(0,\partial H_{\bar{p}}(z_k))^2\geq C_ee_k^{2\theta},
\end{align*}
where \[C_e=\frac{\liminf_{k\to\infty}M_{1,k}}{4c^2(1-\theta)^2\max\bigg\{\left(\frac{L_{\nabla h}}{\underline{\lambda}}+L_{\nabla g}+6\bar{p}\right)^2,  \left(\frac{L_{\nabla h}L_{\nabla g}\overline{\lambda}+1}{\underline{\lambda}}\right)^2\bigg\}}.\]
Applying Lemma~\ref{lem: general convergence rate} completes the proof.\endproof

We now develop sequential convergence rate based on Theorem~\ref{thm: function value rate}.
\begin{theorem}[sequential convergence rate]\label{thm: sequential rate} Let $(x_k)_{k\in\N^*}$ be a sequence generated by BiFRB and define $z_k=(x_{k+1},x_k)$ for $k\in\N^*$. Suppose that all assumptions in Theorem~\ref{thm: finite length} are satisfied and let $x^*$ be the limit given in Theorem~\ref{thm: finite length}(ii). Suppose further that the dominating BiFRB merit function $H_{\bar{p}}$ (recall Definition~\ref{def: merit}) has the generalized concave KL property at $z^*=(x^*,x^*)$ with KL exponent $\theta\in[0,1)$.  Then the following hold:
	
(i) If $\theta=0$, then $z_k$ converges $z^*$ in finite steps.

(ii) If $\theta\in\left(0,1/2\right]$, then there exist $c_1>0$ and $Q_1\in[0,1)$ such that for $k$ sufficiently large
\[\norm{z_k-z^*}\leq c_1Q_1^k.\]

(iii) If $\theta\in\left(1/2,1\right)$, then there exist $c_2>0$ such that for $k$ sufficiently large
\[\norm{z_k-z^*}\leq c_2k^{-\frac{1-\theta}{2\theta-1} }.\]
\end{theorem}
\proof For simplicity, define for $k\in\N^*$
\[e_k=H_{p_k}(z_k)-F(x^*),\delta_k=\sum_{i=k}^\infty\norm{z_{i+1}-z_i}.\]
Then $e_k\to0$ and $\norm{z_k-z^*}\leq \delta_k$. Without loss of generality (recall Lemma~\ref{lem: merit function descent}), assume that $e_k\in[0,1)$ and $e_{k+1}\leq e_k$ for $k\in\N^*$. To obtain the desired results, it suffices to estimate~$\delta_k$.

Assume without loss of generality that $M_{1,k}\geq \liminf_{k\to\infty}M_{1,k}/2$. Invoking Lemma~\ref{lem: merit function descent}(i) and Theorem~\ref{thm: function value convergence}(i) yields that
$\left(\liminf_{k\to\infty}M_{1,k}/2\right)\norm{z_k-z_{k-1}}^2\leq M_{1,k}\norm{z_k-z_{k-1}} \leq H_{p_{k-1}}(z_{k-1})-H_{p_k}(z_k)\leq H_{p_{k-1}}(z_{k-1})-F(x^*)=e_{k-1}$, which in turn implies that
\begin{equation}\label{formula: compare e_k}
\norm{z_k-z_{k-1}}\leq d\sqrt{e_{k-1}},
\end{equation}
where $d=\sqrt{2/\liminf_{k\to\infty}M_{1,k}}$. Moreover, recall from Theorem~\ref{thm: finite length}(i) that there exists $k_0$ such that for $k\geq k_0+1$
\begin{equation}\label{xcv}
\delta_k\leq\norm{z_k-z_{k-1}}+C\tilde{\varphi}(e_k)\leq \norm{z_k-z_{k-1}}+C\tilde{\varphi}(e_{k-1}),
\end{equation}
where the last inequality holds because $e_k\leq e_{k-1}$. Assume without loss of generality that $H_{\bar{p}}$ is associated with desingularizing $\varphi(t)=C^{-1}t^{1-\theta}$ due to the KL exponent assumption. Thus, combined with~(\ref{formula: compare e_k}), inequality~(\ref{xcv}) yields
\begin{equation}\label{formula: sequential rate key inequality}
	\delta_k\leq d\sqrt{e_{k-1}}+e_{k-1}^{1-\theta}.
\end{equation}

Case 1 ($\theta=0$): Appealing to Theorem~\ref{thm: function value rate}, the sequence $e_k$ converges to $0$ in finite steps. Thus~(\ref{formula: sequential rate key inequality}) implies $\delta_k\to0$ in finite steps, so does $(z_k)_{k\in\N^*}$.

Case 2 ($0<\theta\leq\frac{1}{2}$): Clearly $\frac{1}{2}\leq1-\theta<1$, therefore $\sqrt{e_{k-1}}\geq e_{k-1}^{1-\theta}$. Inequality~(\ref{formula: sequential rate key inequality}) and Theorem~\ref{thm: function value rate} imply that for $k$ sufficiently large
 \[\delta_k\leq d\sqrt{e_{k-1}}+\sqrt{e_{k-1}}=(1+d)\sqrt{e_{k-1}}\leq (1+d)\sqrt{\hat{c}_1\hat{Q}_1^{k-1}}.\]

Case 3 ($\frac{1}{2}<\theta<1$): Note that $\sqrt{e_{k-1}}\leq e_{k-1}^{1-\theta}$. Then appealing to~(\ref{formula: sequential rate key inequality}) and Theorem~\ref{thm: function value rate} yields for $k$ sufficiently large
\[\delta_k\leq de_{k-1}^{1-\theta}+e_{k-1}^{1-\theta}=(1+d)e_{k-1}^{1-\theta}\leq (1+d)(k-1)^{-\frac{1-\theta}{2\theta-1}},\]
from which the desired result readily follows.\endproof

\begin{example}[nonconvex feasibility problem] \label{ex: feasibility convergence} Let $C\subseteq\Rn$ be a nonempty, closed and convex set and let $D\subseteq\Rn$ be nonempty and closed. Suppose that $C\cap D\neq\emptyset$ and either $C$ or $D$ is compact. Assume further that both $C$ and $D$ are  semialgebraic.  Consider the minimization problem~(\ref{optimization prblem}) with $f=\delta_D$ and $g=\dist^2(\cdot,C)/2$. Let $(x_k)_{k\in\N^*}$ be a sequence generated by BiFRB. Suppose that Assumption~\ref{assmption: p_k and M_1,k} holds. Then the following hold:
	
	(i) There exists $x^*\in\Rn$ such that $x_k\to x^*$ and $0\in\partial F(x^*)$.
	
	(ii) Suppose additionally that
	\begin{equation}\label{CQ}
		N_C(\proj_C(x^*) )\cap\big(- N_D(x^*) \big)=\{0\}.
	\end{equation}
	Then $x^*\in C\cap D$. Moreover, there exist $Q_1\in(0,1)$ and $c_1>0$ such that \[\norm{x_k-x^*}\leq c_1Q_1^{k},\]
	for $k$ sufficiently large.
\end{example}
\proof (i) By the compactness assumption, the function $f+g$ is coercive. Hence Theorem~\ref{thm: bounded sequence} implies that $(x_k)_{k\in\N^*}$ is bounded. We assume that sets $C,D$ are semialgebraic, then so are functions $f$ and $g$; see~\cite[Section 4.3]{Attouch2010} and~\cite[Lemma 2.3]{attouch2013convergence}.  Then apply Corollary~\ref{cor: convergence on semialgebraic functions}.

(ii) Taking the fact that $0\in\partial F(x^*)$ into account and applying the subdifferential sum rule~\cite[Exercise 8.8]{rockwets}, one gets
\[x^*-\proj_C(x^*)\in-N_D(x^*) .\]
The constraint qualification~(\ref{CQ}) then implies that $x^*=\proj_C(x^*)$ and consequently $x^*\in C$. The set $D$ is assumed to be closed, thus $x^*\in D$. Then a direct application of~\cite[Theorem 5]{chen2020difference} guarantees that $F$ has KL exponent $\theta=1/2$ at $x^*$. The desired result follows immediately from Theorem~\ref{thm: sequential rate}. \endproof

\section{Bregman proximal subproblem formulae}\label{sec:subproblem formulae}
Set $\alpha\geq0,\beta>0$. Fix $\omega\in\Rn$ and $\lambda>0$. In the remainder of this section, let $$h(x)=\alpha\sqrt{1+\norm{x}^2}+\frac{\beta}{2}\norm{x}^2.$$
Define for $u\in\Rn$, 	$p_\lambda(u)=\lambda \omega-\nabla h(u)$ and
\begin{align}
	T_\lambda(u)&=\argmin_{x\in\Rn}\bigg\{f(x)+\ip{x-u}{\omega}+\frac{1}{\lambda}D_h(x,u) \bigg\}\label{formula: Tlam}\\
	&=\argmin_{x\in\Rn}\bigg\{f(x) +\ip{x}{\omega}+\frac{1}{\lambda}h(x)-\frac{1}{\lambda}\ip{\nabla h(u)}{x}\bigg\} \nonumber\nonumber\\
	&=\argmin_{x\in\Rn}\bigg\{\lambda f(x)+\ip{x}{p_\lambda(u)}+h(x)\bigg\}\label{formula: simple Tlambda}.
\end{align}
Clearly Assumption~\ref{assumption: well-defined}(ii) holds with $\sigma=\beta$ and $L_{\nabla h}=\alpha+\beta$ and the Bregman proximal subproblem~(\ref{Algorithm: x_{k+1}}) corresponds to~(\ref{formula: Tlam}) with $\omega=\alpha_k(x_{k-1}-x_k)/\lambda_k+\nabla g(x_k)$, $u=y_k$ and $\lambda=\lambda_{k}$. Note that~(\ref{formula: Tlam}) covers subproblems in~\cite[Algorithm 3.1]{boct2016inertial} and~\cite[Agorithm 1]{boct2016inertiala} as well.

In this section, we develop formulae of $T_\lambda(u)$ with the kernel $h$ chosen as above, supplementing not only Algorithm~\ref{BiFRB} but also~\cite[Algorithm 3.1]{boct2016inertial} and~\cite[Algorithm 1]{boct2016inertiala}.

\subsection{Formula of $l_0$-constrained problems}
Let $D=\{x\in\Rn: \norm{x}_0\leq r,\norm{x}\leq R\}$ for integer an $r\leq n$ and real number $R>0$, where~$(\forall x\in\Rn)~\norm{x}_0=\sum_{i=1}^n|\sgn(x_i)|$ with $\sgn(0)=0$. In this subsection, consider minimization problem
\[\min_{x\in D} g(x),\]
which is~(\ref{optimization prblem}) with $f=\delta_D$. Recall that the Hard-threshold of $x\in\Rn$ with parameter $r$, denoted by $H_r(x)$, is given by
\begin{equation}
(\forall i\in\{1,\ldots,n\})~\big(H_r(x)\big)_i=\begin{cases}x_i, &\text{ if $|x_i|$ belongs to the $r$ largest among $|x_1|,\ldots,|x_n|$},\\
		0,&\text{ otherwise}.\end{cases}
\end{equation}
Immediately, $H_r(cx)=cH_r(x)$ for $c\neq0$. 
The following lemma will be instrumental.
\begin{lemma}\emph{\cite[Proposition 4.3]{luss2013conditional}}\label{lem: inner min solution } Given $0\neq p\in\Rn$ and a positive integer $r<n$, we have
	 \begin{equation*}
		\max_{x\in\Rn}\big\{\ip{p}{x}: \norm{x}=1,~\norm{x}_0\leq r \big\}=\norm{H_r(p)},
	\end{equation*}
with optimal value attained at $x^*=\norm{H_r(p)}^{-1}H_r(p)$.
\end{lemma}

Following a similar route in~\cite[Proposition 5.2]{bolte2018first}, we now prove the desired formula.

\begin{proposition}[subproblem formula]\label{prop: sparse constraints} Let $f=\delta_D$, where $D=\{x\in\Rn: \norm{x}_0\leq r,\norm{x}\leq R\}$ for integer an $r\leq n$ and real number $R>0$. Define
		\begin{equation}\label{formula: Tambda solution with l0 and l2 ball constraints}
		x^*=\begin{cases}
			0,&\text{if }\nabla h(u)=\lambda \omega,\\
			-t^*\norm{H_r\big(p_\lambda(u)\big)}^{-1}H_r\big(p_\lambda(u)\big),&\text{if }\nabla h(u)\neq\lambda \omega,
		\end{cases}
	\end{equation}
where $t^*$ is the unique solution of
\begin{equation}\label{formula: 1d problem}
	\min_{0\leq t\leq R}\bigg\{\varphi(t)=\alpha\sqrt{1+t^2}+\frac{\beta}{2}t^2-\norm{H_r(p_\lambda(u)}t\bigg\},
\end{equation}
which satisfies $t^*=R$ if $\norm{H_r(p_\lambda(u))}\geq\alpha(1+R^2)^{-1/2}+\beta R$; otherwise $t^*\in(0,R)$ is the unique solution of $$\alpha(1+t^2)^{-1/2}t+\beta t-\norm{H_r(p_\lambda(u))}=0.$$
 Then $x^*\in T_\lambda(u)$.
\end{proposition}
\proof Set $p=p_\lambda(u)$ for simplicity. If $\nabla h(u)=\lambda \omega$, then $p=0$ and clearly  $T_\lambda(u)=\argmin_{x\in\Rn}\{h(x): \norm{x}_0\leq r,\norm{x}\leq R\}=0$. Now we consider $\nabla h(u)\neq \lambda \omega$, in which case $p\neq 0$.  Invoking Lemma~\ref{lem: inner min solution }, one has for $t>0$
\begin{align}\label{xcx}
\min_{x\in\Rn}\{\ip{x}{p}: \norm{x}_0\leq r, \norm{x}=t\}&=-\max_{x\in\Rn}\{\ip{x/t}{-tp}: \norm{x/t}\leq r,\norm{x/t}=1\}\nonumber\\
&=-\norm{H_r(-tp)}=-t\norm{H_r(p)},
\end{align}
attained at
\begin{equation}\label{formula: inner min solution}
	\hat{x}=\hat{x}(t)=t\norm{H_r(-tp)}^{-1}H_r(-tp)=-t\norm{H_r(p)}^{-1}H_r(p).
\end{equation}
It is easy to see that optimal value~(\ref{xcx}) and solution~(\ref{formula: inner min solution}) still hold when $t=0$.

Next we show  that $x^*\in T_\lambda (u)$. By~(\ref{formula: simple Tlambda}), $T_\lambda(u)$ is the solution set to
\begin{equation}\label{formula: original optimization problem}
	\min_{x\in\Rn}\bigg\{\ip{x}{p}+\alpha\sqrt{1+\norm{x}^2}+\frac{\beta}{2}\norm{x}^2: \norm{x}_0\leq r, \norm{x}\leq R\bigg\}.
\end{equation}
Consider an arbitrary feasible $x$ of~(\ref{formula: original optimization problem}) and  let $t_x=\norm{x}$. Then $0\leq t_x\leq R$ and
\begin{align*}
\ip{p}{x^*}+\alpha\sqrt{1+\norm{x^*}^2}+\frac{\beta}{2}\norm{x^*}^2&=-t^*\norm{H_r(p)}+\alpha\sqrt{1+(t^*)^2}+\frac{\beta}{2}(t^*)^2\\
&=\varphi(t^*)\leq\varphi(t_x)=-\norm{H_r(p)}t_x+\alpha\sqrt{1+t_x^2}+\frac{\beta}{2}t_x^2\\
&\leq\ip{p}{x}+ \alpha\sqrt{1+\norm{x}^2}+\frac{\beta}{2}\norm{x}^2,
\end{align*}
where the first and last inequalities hold due to~(\ref{formula: 1d problem}) and (\ref{xcx}) respectively.

Finally, we turn to~(\ref{formula: 1d problem}). Noticing the strong convexity of $\varphi$ and applying optimality condition, one gets \[-\alpha(1+t^2)^{-1/2}t-\beta t+\norm{H_r(p)}=-\varphi^\prime(t^*)\in N_{[0,R]}(t^*),\]
So $$t^*=R\Leftrightarrow -\varphi^\prime(R)\in N_{[0,R]}(R)=[0,\infty)\Leftrightarrow\varphi^\prime(R)\leq 0\Leftrightarrow\norm{H_r(p)}\geq\alpha(1+R^2)^{-1/2}+\beta R,$$ and  $$t^*=0\Leftrightarrow-\varphi^\prime(0)\in N_{[0,R]}(0)=(-\infty,0]\Leftrightarrow\norm{H_r(p)}\leq0\Leftrightarrow H_r(p)=0,$$
which never occurs when $p\neq0$.
If $0<\norm{H_r(p)}<\alpha(1+R^2)^{-1/2}+\beta R$, then $\varphi^\prime(0)<0$ and $\varphi^\prime(R)>0$, meaning that there exists unique $t^*\in(0,R)$ such that $\varphi^\prime(t^*)=0$ due to the strong convexity of $\varphi$. \endproof
\begin{remark}\label{remark: root-finding} When $0<\norm{H_r(p_\lambda(u)}<\alpha(1+R^2)^{-1/2}+\beta R$, $t^*$ in Proposition~\ref{prop: sparse constraints} is indeed one of the four roots of quartic equation
	\[\beta^2t^4-2\beta\norm{H_r(p_\lambda(u))}t^3+\left( \norm{H_r(p_\lambda(u))}^2+\beta^2-\alpha^2\right)t^2-2\beta\norm{H_r(p_\lambda(u))}t+\norm{H_r(p_\lambda(u))}^2=0.\]
Although quartic equations admit closed-form solutions, it is virtually impossible to determine beforehand the branch of solutions to which $t^*$ belongs, due to parameters $\alpha$, $\beta$ and $\norm{H_r(p_\lambda(u))}$.
	To find $t^*$ efficiently, one may employ numerical root-finding algorithms. For instance, apply bisection method to identify a subinterval of $[0,R]$ in which the desired root $t^*$ lies, then appeal to Newton's method; see, e.g.,~\cite[Chapter 2]{BFB10E}.
\end{remark}

\subsection{Formula of problems penalized by a positively homogeneous convex function}

Recall that $f:\Rn\to\overline{\R}$ is positively homogeneous if $f(tx)=tf(x)$ for every $x\in\Rn$ and $t>0$. In this subsection, we provide subproblem formula for problems of the form
\[\min_{x\in\Rn} f(x)+g(x), \]
where $f:\Rn\to\overline{\R}$ is convex and positively homogeneous with known proximal mapping, which is a special case of problem~(\ref{optimization prblem}).

The result below generalizes the technique of Bolte-Sabach-Teboulle-Vaisbourd~\cite[Proposition 5.1]{bolte2018first}. However one cannot recover their result as our kernel is different.

\begin{proposition}\label{prop: formula for p.h. penalized problems} Let $f:\Rn\to\overline{\R}$ be proper, lsc and convex. Suppose that $f$ is positively homogeneous. Then the following hold:

(i) $(\forall x\in\dom\partial f)$ $(\forall t>0)$ $\partial f(x)=\partial f(tx)$.

(ii) \emph{(subproblem formula)} $(\forall u\in\Rn)$ $T_\lambda(u)=t^*\prox_{\lambda f}(-p_\lambda(u))$, where $t^*$ is the unique root of the equation
\[1-\alpha t\big(1+t^2\norm{\prox_{\lambda f}(-p_\lambda(u))} \big)^{-1/2}-\beta t=0. \]

\end{proposition}
\proof (i) $v\in\partial f(x)\Leftrightarrow (\forall y\in\Rn)~f(y)\geq f(x)+\ip{v}{y-x}\Leftrightarrow (\forall y\in\Rn)~(\forall t>0)~f(ty)\geq f(tx)+\ip{v}{ty-tx}\Leftrightarrow v\in\partial f(tx)$.

(ii) For simplicity, let $p=p_\lambda(u)$. Note that $T_\lambda(p)$ is single-valued by strong convexity. Then invoking optimality condition and subdifferential sum rule yields
\begin{equation*}
x^*=T_\lambda(u)\Leftrightarrow 0\in\lambda\partial f(x^*)+p+\nabla h(x)=\lambda \partial f(x^*)+p+[\alpha(1+\norm{x^*}^2)^{-1/2}+\beta]x^*.
\end{equation*}
Define $t^*=[\alpha(1+\norm{x^*}^2)^{-1/2}+\beta]^{-1}$ and $v=x^*/t^*$. Then $t^*>0$ and the above inclusion together with statement(i) implies that
\begin{align*}
	0&\in\lambda \partial f(x^*)+p+v=\lambda\partial f(t^*v)+p+v=\lambda \partial f(v)+v+p\\
	\Leftrightarrow -p&\in(\Id+\lambda\partial f)(v)\Leftrightarrow v=(\Id+\lambda\partial f)^{-1}(-p)=\prox_{\lambda f}(-p).
\end{align*}
In turn, $t^*$ satisfies
\[v[1-\alpha t^*(1+(t^*)^2\norm{v}^2)^{-1/2}-\beta t^*]=v-[\alpha(1+(t^*)^2\norm{v}^2)^{-1/2}+\beta]t^*v=v-\frac{x^*}{t^*}=v-v=0.\]
If $v\neq0$, then clearly $1-\alpha t^*(1+(t^*)^2\norm{v}^2)^{-1/2}-\beta t^*=0$ as claimed; if $v=0$, then $x^*=0$ and $t^*=1/(\alpha+\beta)$, which is consistent with statement(ii). To see the uniqueness, let $\varphi(t)=\alpha(1+\norm{v}^2t^2)^{1/2}/\norm{v}^2+\beta t^2/2-t$ for $t\in\R$, which is strongly convex. Simple calculus shows that $t^*$ is stationary point of $\varphi$, thus unique. \endproof
\begin{remark}  Similarly to Proposition~\ref{prop: sparse constraints} and Remark~\ref{remark: root-finding}, one may find $t^*$ efficiently by appealing to root-finding algorithms.
\end{remark}

Equipped with Proposition~\ref{prop: formula for p.h. penalized problems}, we now turn to concrete problems. Recall that the Soft-threshold of $x\in\Rn$ with parameter $\lambda$ is $S_\lambda(x)=\prox_{\lambda \norm{\cdot}_1}(x)$, which satisfies $$(\forall i\in\{1,\ldots,n\})~(S_\lambda(x))_i=\max(|x_i|-\lambda,0)\sgn(x_i),$$
see, e.g.,~\cite[Example 6.8]{beck2017first}.
\begin{corollary}[$l_1$-penalized problem] Let $f=\norm{\cdot}_1$. Then $	T_\lambda(u)=-t^*S_\lambda\left(p_\lambda(u)\right)$,
where $t^*$ is the unique solution of $1-\alpha t\left(1+t^2\norm{S_\lambda\left(p_\lambda(u)\right)}^2\right)^{-1/2}-\beta t=0$.
\end{corollary}

\begin{corollary}[$l_\infty$-penalized problem] Let $f=\norm{\cdot}_\infty$. Then $T_\lambda(u)=t^*v$, where $$v=p_\lambda(u)-\lambda\proj\left(p_\lambda(u)/\lambda;\Ball_{\norm{\cdot}_1}\right),$$ and $t^*$ is the unique solution of $1-\alpha t(1+t^2\norm{v}^2)^{-1/2}-\beta t=0$.
\end{corollary}
\begin{remark} The projection $\proj(\cdot;\Ball_{\norm{\cdot}_1})$ admits a closed-form formula; see, e.g.,~\cite[Example 6.33]{beck2017first}.
\end{remark}

\section{Numerical experiments}\label{sec:numerical}
In the remainder, let $A\in\R^{m\times n}$, $b\in\R^m$, $r\in\N\backslash\{0\}$, and $R>0$. Define $C=\{x\in\Rn: Ax=b\}$ and $D=\{x\in\Rn: \norm{x}_0\leq r,\norm{x}\leq R\}$. Clearly $C$ is semialgebraic. The set $D$ can be represented as
\[D=\{x\in\Rn: \norm{x}_0\leq r \}\cap\{x\in\Rn: \norm{x}^2\leq R^2\}, \]
which is an intersection of semialgebraic sets; see, e.g,~\cite[Formula 27(d)]{bauschke2014restricted}, thus semialgebraic. Consider
\begin{equation}\label{problem: nonconvex fea}
	\min_{x\in D}\frac{1}{2}\dist^2(x,C),
\end{equation}
which corresponds to~(\ref{optimization prblem}) with $f=\delta_D$ and $g=\dist^2(\cdot,C)/2$. In turn, Example~\ref{ex: feasibility convergence} ensures that BiFRB and its Euclidean variant iFRB converge to a stationary point of~(\ref{problem: nonconvex fea}).

We shall benchmark BiFRB and iFRB against the Douglas-Rachford~(DR)~\cite{li2016douglas}, Forward-reflected-backward~(FRB)~\cite{wang2021malitsky,malitsky2020FRB} and inertial Tseng's~(iTseng)~\cite{boct2016inertial} methods with $R=1,1000$ and $r=\lceil m/5\rceil$. These splitting methods are
known to converge globally to a stationary point of~(\ref{optimization prblem}) under the same assumptions on $F$; see~\cite[Theorems 1--2, Remark 4, Corollary 1]{li2016douglas},~\cite[Theorem 3.9]{wang2021malitsky} and~\cite[Theorem 3.1]{boct2016inertial}, respectively.

BiFRB, iFRB and FRB are implemented with stepsize rules given by Proposition~\ref{prop: non-Euclidean fixed stepsize}, Corollary~\ref{cor: Euclidean fixed stepsize} and~\cite[Theorem 3.9]{wang2021malitsky} respectively, and terminated when
\begin{equation}\label{formula: termination condition}
	\frac{\max\big\{\norm{x_{k+1}-x_k },\norm{x_k-x_{k-1}} \big\} }{\max\big\{1,\norm{x_k},\norm{x_{k-1}}\big\}}<10^{-10}.
\end{equation}
Moreover, BiFRB uses the kernel $h(x)=0.1\cdot \sqrt{1+\norm{x}^2}+2.51\cdot\norm{x}^2/2$.
We apply DR with stepsize and termination condition given by~\cite[Section 5]{li2016douglas} with the same tolerance $10^{-10}$; while iTseng employs a stepsize given by~\cite[Lemma 3.3]{boct2016inertial} and termination condition~(\ref{formula: termination condition}). As for inertial parameter, we use $\alpha_1=0.9$ for BiFRB and $\alpha_2=0.49$ for both of iFRB and iTseng. Finally, each algorithm is initialized at the origin and is equipped with the stepsize heuristic described in~\cite[Remark 4]{li2016douglas}.

Simulation results are presented in Table~\ref{num. results} below. Our problem data is generated through creating random matrices $A\in\R^{m\times n}$ with entries following the standard Gaussian distribution. For each problem of the size~$(m,n)$, we randomly generate 50 instances and report ceilings of the average number of iterations (iter) and the minimal objective value at termination ($\text{fval}_{\text{min}}$). As problem~(\ref{problem: nonconvex fea}) has optimal value zero, $\text{fval}_{\text{min}}$ indicates whether an algorithm actually hits a global minimizer rather than simply converging to a stationary point. We observed the following:
\begin{itemize}
	\item[-] BiFRB is the most advantageous method on ``bad" problems (small $m$), but it is not robust.
	\item[-] DR tends to have the smallest function value at termination, however, it can converge very slowly.
	\item[-] In most cases, iFRB has the smallest number of iterations with fair function values at termination.
\end{itemize}

\section{Conclusions}\label{sec:conclusion and future work}
We proved convergence of BiFRB (Algorithm~\ref{BiFRB}) by using the generalized concave KL property and a general framework for decreasing property. The resulting stepsize condition is independent of inertial parameter, which is less restrictive compared to other algorithms with the same setting. In turn, a question of Malitsky and Tam~\cite[Section 7]{malitsky2020FRB} is answered. We believe that our approach could be useful for other splitting algorithms in designing more implementable stepsize rules in the absence of convexity. We also conducted convergence rate analysis on both function value and the actual sequence. Formulae for our Bregman subproblems are provided as well, which supplements not only BiFRB but also~\cite{boct2016inertial,boct2016inertiala}. To end this paper, let us now present several directions
for future work:
\begin{itemize}
	\item[-] It is tempting to see whether the global Lipschitz assumption on gradients $\nabla g$ and $\nabla h$ in Assumption~\ref{Standing assumption} can be relaxed in a fashion similar to~\cite{bolte2018first}.
	\item[-] Another interesting question is to explore the reason why BiFRB demonstrates inconsistent performance when problem size varies; see Section~\ref{sec:numerical}.
	\item[-] Now that the range of admissible inertial parameter is extended (recall Proposition~\ref{prop: non-Euclidean fixed stepsize}), it would be interesting to see if there exists an optimal accelerating scheme.
\end{itemize}

\begin{landscape}
\begin{table}\caption{Comparing BiFRB and iFRB to FRB, DR and iTseng.\label{num. results}}
\begin{center}
	\begin{tabular}{|c c|c c|c c|c c|c c|c c|}
	\hline
	Size & $(R=1)$&  BiFRB  &&  iFRB&& FRB && DR&& iTseng&\\
	\hline
	m&n&iter&$\text{fval}_{\text{min}}$&iter&$\text{fval}_{\text{min}}$&iter&$\text{fval}_{\text{min}}$&iter&$\text{fval}_{\text{min}}$&iter&$\text{fval}_{\text{min}}$\\
	\hline
	100 &4000&50 &0.03251&93  &0.03251
	&631  &0.03929
	&860  &0.02819
	&1367  &0.03149 \\
	\hline
	100&5000  &52  &0.02413
	&115  &0.02413
	&733  &0.02875
	&684  &0.02192
	&1632  &0.02268\\
	\hline
	100&6000  &57  &0.01369
	&145  &0.01369
	&911  &0.01794
	&683  &0.01253
	&1970  &0.0133\\
	\hline
	200&4000  &870  &0.2923
	&39  &0.2745
	&190  &0.2769
	&5466  &0.2711
	&448  &0.2746\\
	\hline
	200&5000  &1198  &0.2367
	&41  &0.2095
	&231  &0.2157
	&5864  &0.2073
	&546  &0.2097\\
	\hline
	200&6000  &81  &0.2306
	&43  &0.2259
	&257  &0.2292
	&4199  &0.2236
	&620  &0.2264\\
	\hline
	300&4000  &885  &1.051
	&32  &1.036
	&105  &1.038
	&5497  &1.036
	&253  &1.048\\
	\hline
	300&5000  &925  &0.7481
	&34  &0.7044
	&124  &0.7071
	&5673  &0.7032
	&302  &0.7051\\
	\hline
	300 &6000  &1407  &0.583
	&36  &0.5546
	&144  &0.555
	&6173  &0.5541
	&352  &0.5566\\
	\hline
	\hline
	Size & $(R=1000)$&  BiFRB  &&  iFRB&& FRB && DR&& iTseng&\\
	\hline
	m&n&iter&$\text{fval}_{\text{min}}$&iter&$\text{fval}_{\text{min}}$&iter&$\text{fval}_{\text{min}}$&iter&$\text{fval}_{\text{min}}$&iter&$\text{fval}_{\text{min}}$\\
	\hline
	$100$& $4000$ &1873 &0.006609
	&1210 &0.00365
	&7376 &0.00816
	&2194 &4e-21
	&10001 &0.01671\\
	\hline
	$100$& $5000$ &574 &0.00278
	&1427 &0.002244
	&8675 &1.098e-05
	&2919 &1.626e-20
	&10001 &0.0158\\
	\hline
	$100$ & $6000$ &790 &0.003827
	&1696 &0.001261
	&9394 &0.0002952
	&2344 &3.438e-20
	&10001 &0.01066\\
	\hline
	$200$ & $4000$  &5842  &2.992
	&750  &2.442e-18
	&8223  &0.007688
	&1038  &1.106e-20
	&9797  &0.1094\\
	\hline
	$200$ &$5000$  &5938  &0.6547
	&859  &2.876e-05
	&7257  &0.02581
	&1242  &2.082e-20
	&9997  &0.07125\\
	\hline
	$200$& $6000$  &5666  &0.917
	&1086  &0.0001211
	&7312  &0.006726
	&1487  &9.129e-21
	&10001  &0.07524 \\
	\hline
	$300$&$4000$  &7164  &3.322
	&619  &2.196e-18
	&4620  &1.462e-05
	&753  &9.795e-21
	&7722  &0.2417\\
	\hline
	$300$&$5000$  &5267  &2.299
	&693  &2.175e-18
	&6352  &0.02936
	&880  &5.446e-20
	&9194  &0.1702 \\
	\hline
	$300$&$6000$  &5724  &5.06
	&750  &4.71e-18
	&8644  &0.1109
	&1027  &2.175e-20
	&9823  &0.1423 \\
	\hline
\end{tabular}
	
\end{center}
\end{table}
\end{landscape}
\section*{Acknowledgments}
XW and ZW were partially supported by NSERC Discovery Grants.

\bibliographystyle{siam}
\bibliography{nonconvexFRB}

\end{document}